Sachin Malik

Neeraj Kumar

Florentin Smarandache

*(Editors)*

USES OF SAMPLING

TECHNIQUES

&

INVENTORY CONTROL

WITH CAPACITY

CONSTRAINTS

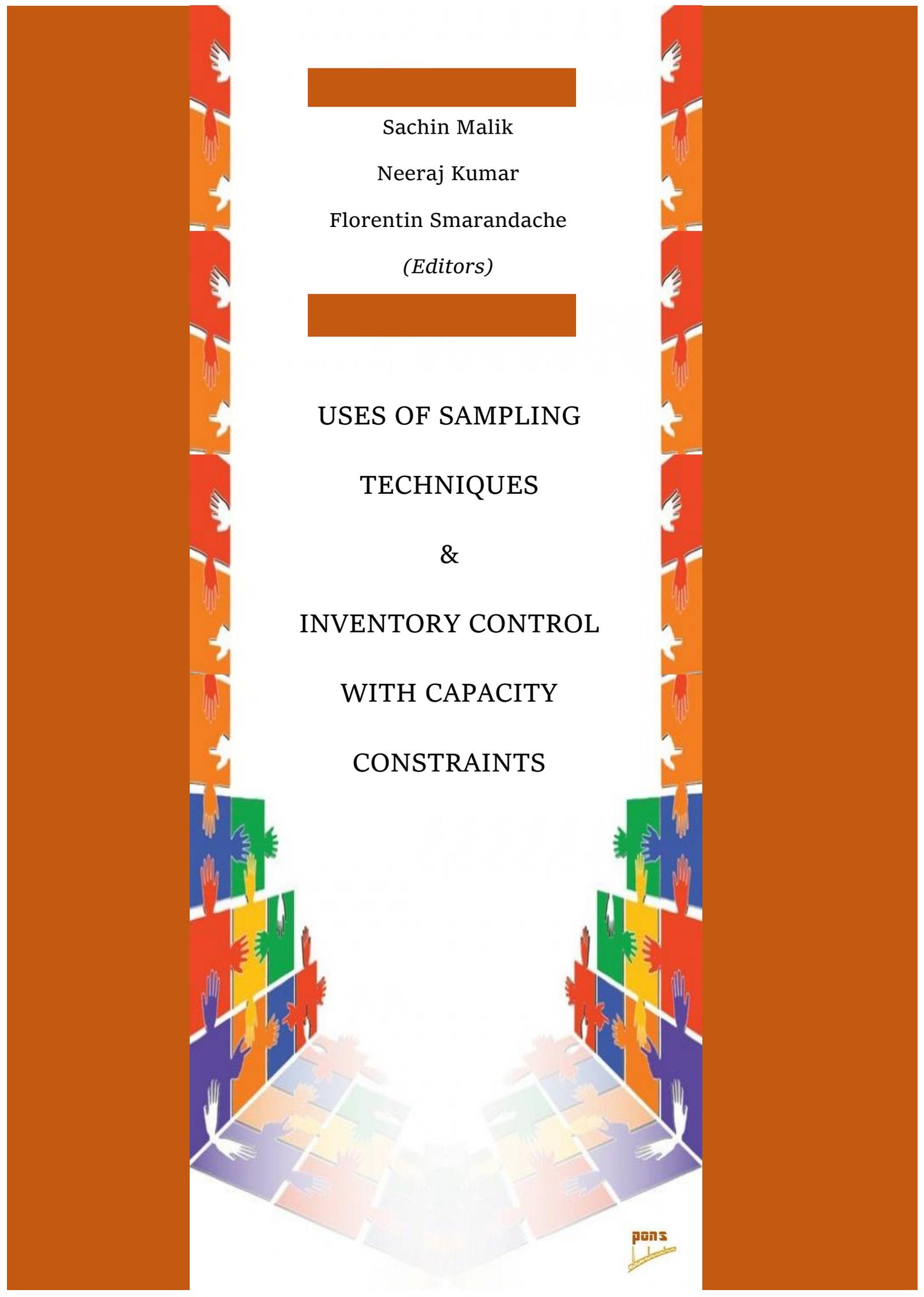

Sachin Malik, Neeraj Kumar, Florentin Smarandache

*Editors*

USES OF SAMPLING TECHNIQUES
&
INVENTORY CONTROL WITH CAPACITY
CONSTRAINTS

*Peer-Reviewers:*

**Nassim Abbas**, Communicating and Intelligent System Engineering Laboratory, Faculty of Electronics and Computer Science University of Science and Technology Houari Boumediene 32, El Alia, Bab Ezzouar, 16111, Algiers, Algeria.

**Mumtaz Ali**, Department of Mathematics, Quaid-i-Azam University, Islamabad, 44000, Pakistan.

**Said Broumi**, Faculty of letters and Humanities, Hay El Baraka Ben M'sik Casablanca, B.P. 7951, University of Hassan II, Casablanca, Morocco.

**Jun Ye**, Department of Electrical and Information Engineering, Shaoxing University, No. 508 Huancheng West Road, Shaoxing, Zhejiang Province 312000, P. R. China.

# USES OF SAMPLING TECHNIQUES
# &
# INVENTORY CONTROL WITH CAPACITY CONSTRAINTS

*Editors:*

**Sachin Malik, Neeraj Kumar**
Department of Mathematics, SRM University
Delhi NCR, Sonepat - 131029, India

**Florentin Smarandache**
Department of Mathematics, University of New Mexico,
Gallup, NM 87301, USA

*Pons Editions*

Brussels, 2016

# Contents

# *Preface*

The main aim of the present book is to suggest some improved estimators using auxiliary and attribute information in case of simple random sampling and stratified random sampling and some inventory models related to capacity constraints.

This volume is a collection of five papers, written by six co-authors (listed in the order of the papers): Dr. Rajesh Singh, Dr. Sachin Malik, Dr. Florentin Smarandache, Dr. Neeraj Kumar, Mr. Sanjey Kumar & Pallavi Agarwal.

In the first chapter authors suggest an estimator using two auxiliary variables in stratified random sampling for estimating population mean. In second chapter they proposed a family of estimators for estimating population means using known value of some population parameters. In Chapter third an almost unbiased estimator using known value of some population parameter(s) with known population proportion of an auxiliary variable has been used. In Chapter four the authors investigates a fuzzy economic order quantity model for two storage facility. The demand, holding cost, ordering cost, storage capacity of the own - warehouse are taken as trapezoidal fuzzy numbers. And in Chapter five a two-warehouse inventory model deals with deteriorating items, with stock dependent demand rate and model affected by inflation under the pattern of time value of money over a finite planning horizon. Shortages are allowed and partially backordered depending on the waiting time for the next replenishment. The purpose of this model is to minimize the total inventory cost by using the genetic algorithm.

This book will be helpful for the researchers and students who are working in the field of sampling techniques and inventory control.

# An Improved Suggestion in Stratified Random Sampling Using Two Auxiliary Variables

**Rajesh Singh, Sachin Malik\* and Florentin Smarandache\*\***

\*Department of Mathematics, SRM University

Delhi NCR, Sonepat- 131029, India

Department of Statistics, Banaras Hindu University

Varanasi-221005, India

\*\*Department of Mathematics, University of New Mexico

Gallup, NM 87301, USA

\*Corresponding Author, sachinkurava999@gmail.com

**Abstract**

In this paper, we suggest an estimator using two auxiliary variables in stratified random sampling following Malik and Singh [12]. The propose estimator has an improvement over mean per unit estimator as well as some other considered estimators. Expressions for bias and MSE of the estimator are derived up to first degree of approximation. Moreover, these theoretical findings are supported by a numerical example with original data.

## 1. Introduction

The problem of estimating the population mean in the presence of an auxiliary variable has been widely discussed in finite population sampling literature. Out of many ratio, product and regression methods of estimation are good examples in this context. Diana [2] suggested a class of estimators of the population mean using one auxiliary variable in the stratified random sampling and examined the MSE of the estimators up to the k[th] order of approximation. Kadilar and Cingi [3], Singh et al. [7], Singh and Vishwakarma [8],Koyuncu and Kadilar [4] proposed estimators in stratified random sampling. Singh [9] and Perri [6] suggested some ratio cum product estimators in simple random sampling. Bahl and Tuteja [1] and Singh et al. [11] suggested some exponential ratio type estimators. In this chapter, we suggest some exponential-type estimators using the auxiliary information in the stratified random sampling.

Consider a finite population of size N and is divided into L strata such that $\sum_{h=1}^{L} N_h = N$ where $N_h$ is the size of $h^{th}$ stratum (h=1,2,...,L). We select a sample of size $n_h$ from each stratum by simple random sample without replacement sampling such that $\sum_{h=1}^{L} n_h = n$, where $n_h$ is the stratum sample size. A simple random sample of size $n_h$ is drawn without replacement from the $h^{th}$ stratum such that $\sum_{h=1}^{L} n_h = n$. Let ($y_{hi}$, $x_{hi}$, $z_{hi}$) denote the observed values of y, x, and z on the $i^{th}$ unit of the $h^{th}$ stratum, where i=1, 2, 3...$N_h$.

To obtain the bias and MSE, we write

$$\bar{y}_{st} = \sum_{h=1}^{L} w_h \bar{y}_h = \bar{Y}(1+e_0), \quad \bar{x}_{st} = \sum_{h=1}^{L} w_h \bar{x}_h = \bar{X}(1+e_1), \quad \bar{z}_{st} = \sum_{h=1}^{L} w_h \bar{z}_h = \bar{Z}(1+e_2)$$

Such that,

$$E(e_0) = E(e_0) = E(e_0) = 0$$

$$V_{rst} = \sum_{h=1}^{L} w_h^{r+s+t} \frac{E\left[\left(\bar{y}_h - \bar{Y}\right)^r \left(\bar{x}_h - \bar{X}\right)^s \left(\bar{z}_h - \bar{Z}\right)^t\right]}{\bar{Y}^r \bar{X}^s \bar{Z}^t}$$

where,

$$\bar{y}_{st} = \sum_{h=1}^{L} w_h \bar{y}_h \ , \ \bar{y}_h = \frac{1}{n_h} \sum_{i=1}^{n_h} \bar{y}_{hi} \ , \ \bar{Y}_h = \frac{1}{N_h} \sum_{i=1}^{n_h} \bar{Y}_{hi}$$

$$Y = \bar{Y}_{st} = \sum_{h=1}^{L} w_h \bar{Y}_h, \quad w_h = \frac{N_h}{N}$$

and

$$V(\bar{y}_{st}) = \bar{Y}^2 V_{200} \tag{1.1}$$

Similar expressions for X and Z can also be defined.

And $\quad E\left(e_0^2\right) = \dfrac{\sum\limits_{h=1}^{L} W_h^2 f_h S_{yh}^2}{\overline{Y}^2} = V_{200}$, $\qquad\qquad E\left(e_1^2\right) = \dfrac{\sum\limits_{h=1}^{L} W_h^2 f_h S_{xh}^2}{\overline{X}^2} = V_{020}$,

$E\left(e_2^2\right) = \dfrac{\sum\limits_{h=1}^{L} W_h^2 f_h S_{zh}^2}{\overline{Z}^2} = V_{002}$, $\qquad\qquad E\left(e_0 e_1\right) = \dfrac{\sum\limits_{h=1}^{L} W_h^2 f_h S_{yxh}^2}{\overline{Y}\,\overline{X}} = V_{110}$,

$E\left(e_0 e_2\right) = \dfrac{\sum\limits_{h=1}^{L} W_h^2 f_h S_{yzh}^2}{\overline{Y}\,\overline{Z}} = V_{101}$, $\qquad$ and $E\left(e_1 e_2\right) = \dfrac{\sum\limits_{h=1}^{L} W_h^2 f_h S_{xzh}^2}{\overline{X}\,\overline{Z}} = V_{011}$,

where ,

$S_{yh}^2 = \sum\limits_{i=1}^{N_h} \dfrac{\left(y_h - \overline{Y}_h\right)^2}{N_h - 1}$, $\qquad\qquad S_{xh}^2 = \sum\limits_{i=1}^{N_h} \dfrac{\left(x_h - \overline{X}_h\right)^2}{N_h - 1}$

$S_{zh}^2 = \sum\limits_{i=1}^{N_h} \dfrac{\left(z_h - \overline{Z}_h\right)^2}{N_h - 1}$, $\qquad\qquad S_{yxh} = \sum\limits_{i=1}^{N_h} \dfrac{\left(x_h - \overline{X}_h\right)\left(y_h - \overline{Y}_h\right)}{N_h - 1}$

$S_{yzh} = \sum\limits_{i=1}^{N_h} \dfrac{\left(z_h - \overline{Z}_h\right)\left(y_h - \overline{Y}_h\right)}{N_h - 1}$, $\qquad\qquad S_{xzh} = \sum\limits_{i=1}^{N_h} \dfrac{\left(x_h - \overline{X}_h\right)\left(z_h - \overline{Z}_h\right)}{N_h - 1}$

And,

$f_h = \dfrac{1}{n_h} - \dfrac{1}{N_h}$

## 2. Estimators in literature

In order to have an estimate of the study variable y, assuming the knowledge of the population proportion P, Naik and Gupta [5] and Singh et al. [11] respectively proposed following estimators

$t_1 = \overline{y}_{st}\left(\dfrac{\overline{X}}{\overline{x}_{st}}\right)$

$\qquad\qquad\qquad\qquad\qquad\qquad\qquad\qquad\qquad\qquad$ (2.1)

$t_2 = \overline{y}_{st}\exp\left(\dfrac{\overline{X} - \overline{x}_{st}}{\overline{X} + \overline{x}_{st}}\right)$

$\qquad\qquad\qquad\qquad\qquad\qquad\qquad\qquad\qquad\qquad$ (2.2)

The MSE expressions of these estimators are given as

$$\text{MSE}(t_1) = \overline{Y}^2 \left[ V_{200} + V_{020} - 2V_{110} \right] \tag{2.3}$$

$$\text{MSE}(t_2) = \overline{Y}^2 \left[ V_{200} + \frac{V_{020}}{4} - V_{110} \right] \tag{2.4}$$

When the information on the two auxiliary variables is known, Singh [10] proposed some ratio cum product estimators in simple random sampling to estimate the population mean of the study variable y.

Motivated by Singh [10] and Singh et al. [7], Singh and kumar propose some estimators in stratified sampling as

$$t_3 = \overline{y}_{st} \exp\left[ \frac{\overline{X} - \overline{x}_{st}}{\overline{X} + \overline{x}_{st}} \right] \exp\left[ \frac{\overline{Z} - \overline{z}_{st}}{\overline{Z} + \overline{z}_{st}} \right] \tag{2.5}$$

$$t_4 = \overline{y}_{st} \exp\left[ \frac{\overline{x}_{st} - \overline{X}}{\overline{x}_{st} + \overline{X}} \right] \exp\left[ \frac{\overline{z}_{st} - \overline{Z}}{\overline{z}_{st} + \overline{Z}} \right] \tag{2.6}$$

$$t_5 = \overline{y}_{st} \exp\left[ \frac{\overline{X} - \overline{x}_{st}}{\overline{X} + \overline{x}_{st}} \right] \exp\left[ \frac{\overline{z}_{st} - \overline{Z}}{\overline{z}_{st} + \overline{Z}} \right] \tag{2.7}$$

$$t_6 = \overline{y}_{st} \exp\left[ \frac{\overline{x}_{st} - \overline{X}}{\overline{x}_{st} + \overline{X}} \right] \exp\left[ \frac{\overline{Z} - \overline{z}_{st}}{\overline{Z} + \overline{z}_{st}} \right] \tag{2.8}$$

The MSE equations of these estimators can be written as

$$\text{MSE}(t_3) = \overline{Y}^2 \left[ V_{200} + \frac{V_{020}}{4} + \frac{V_{002}}{4} - V_{110} - V_{101} + \frac{V_{011}}{2} \right] \tag{2.9}$$

$$\text{MSE}(t_4) = \overline{Y}^2 \left[ V_{200} + \frac{V_{020}}{4} + \frac{V_{002}}{4} + V_{110} + V_{101} + \frac{V_{011}}{2} \right] \tag{2.10}$$

$$\text{MSE}(t_5) = \overline{Y}^2 \left[ V_{200} + \frac{V_{020}}{4} + \frac{V_{002}}{4} - V_{110} + V_{101} - \frac{V_{011}}{2} \right] \tag{2.11}$$

$$MSE(t_6) = \overline{Y}^2 \left[ V_{200} + \frac{V_{020}}{4} + \frac{V_{002}}{4} + V_{110} - V_{101} - \frac{V_{011}}{2} \right]$$

(2.12)

When there are two auxiliary variables, the regression estimator of $\overline{Y}$ will be

$$t_7 = \overline{y}_{st} + b_{1h}\left(\overline{X} - \overline{x}_{st}\right) + b_{2h}\left(\overline{Z} - \overline{z}_{st}\right)$$

(2.13)

Where $b_{1h} = \frac{s_{yx}}{s_x^2}$ and $b_{2h} = \frac{s_{yz}}{s_z^2}$. Here $s_x^2$ and $s_z^2$ are the sample variances of x and z respectively,

$s_{yx}$ and $s_{yz}$ are the sample covariance's between y and x and between z respectively. The MSE

expression of this estimator is:

$$MSE(t_7) = \sum_{h=1}^{L} W_h^2 f_h S_{yh}^2 \left(1 - \rho_{yxh}^2 - \rho_{yzh}^2 + 2\rho_{yxh}\rho_{yzh}\rho_{xzh}\right)$$

(2.14)

## 3. The proposed estimator

Following Malik and Singh [12], we propose an estimator using information on two auxiliary

attributes as

$$t_p = \overline{y}_{st} \exp\left[\frac{\overline{X} - \overline{x}_{st}}{\overline{X} + \overline{x}_{st}}\right]^{m1} \exp\left[\frac{\overline{Z} - \overline{z}_{st}}{\overline{Z} + \overline{z}_{st}}\right]^{m2} + b_{1h}\left(\overline{X} - \overline{x}_{st}\right) + b_{2h}\left(\overline{Z} - \overline{z}_{st}\right)$$

(3.1)

Expressing equation (3.1) in terms of e's, we have

$$t_p = \overline{Y}\left(1 + e_0\right)\left\{ \exp\left(\frac{-e_1}{2 + e_1}\right)^{m_1} \exp\left(\frac{-e_2}{2 + e_2}\right)^{m_2} \right\} - b_{1h}\overline{X}e_1 - b_{2h}\overline{Z}e_2$$

$$= \overline{Y}\left[1 + e_0 - \frac{m_1 e_1}{2} + \frac{m_1^2 e_1^2}{4} - \frac{m_2 e_2}{2} - \frac{m_1 m_2 e_1 e_2}{4} + \frac{m_2^2 e_2^2}{4} - \frac{m_2 e_0 e_2}{2} - \frac{m_1 e_0 e_1}{2}\right]$$

$$- b_{1h}e_1\overline{X} - b_{2h}e_2\overline{Z}$$

(3.2)

Squaring both sides of (3.2) and neglecting the term having power greater than two, we have

$$\left(t_p - \overline{Y}\right)^2 = \left\{\overline{Y}\left[e_0 - \frac{m_1 e_1}{2} - \frac{m_2 e_2}{2}\right] - b_{1h} e_1 \overline{X} - b_{2h} e_2 \overline{Z}\right\}^2 \tag{3.3}$$

Taking expectations of both the sides of (3.3), we have the mean squared error of $t_p$ up to the first degree of approximation as

$$MSE(t_p) = \overline{Y}^2\left[V_{200} + P_1\right] + P_2 - \overline{Y}P_3 \tag{3.4}$$

Where,

$$
\left.
\begin{aligned}
P_1 &= \frac{m_1^2 V_{020}}{4} + \frac{m_2^2 V_{002}}{4} + \frac{m_1 m_2 V_{011}}{2} - m_1 V_{110} - m_2 V_{101} \\
P_2 &= B_{1h}^2 V_{020} + B_{2h}^2 V_{002} + 2B_{1h} B_{2h} V_{011} \\
P_3 &= -2B_{1h} V_{110} - 2B_{2h} V_{101} + m_1 B_{1h} V_{020} + m_1 B_{2h} V_{011} + m_2 B_{1h} V_{011} + m_2 B_{2h} V_{002}
\end{aligned}
\right\} \tag{3.5}
$$

Where, $B_{1h} = \dfrac{\sum\limits_{h=1}^{L} W_h^2 f_h \rho_{yxh} S_{yh} S_{xh}}{\sum\limits_{h=1}^{L} W_h^2 f_h S_{xh}^2}$ and $B_{2h} = \dfrac{\sum\limits_{h=1}^{L} W_h^2 f_h \rho_{yzh} S_{yh} S_{zh}}{\sum\limits_{h=1}^{L} W_h^2 f_h S_{zh}^2}$

The optimum values of $m_1$ and $m_2$ will be

$$
\left.
\begin{aligned}
m_1 &= \frac{4\left[B_{1h} V_{011} V_{002} + B_{2h} V_{011}^2 - B_{1h} V_{020} V_{002} - B_{2h} V_{011} V_{002}\right]}{\overline{Y}\left[V_{020} V_{002} - V_{011}^2\right]} \\
m_2 &= \frac{4\left[B_{1h} V_{011} V_{020} + B_{2h} V_{011}^2 - B_{1h} V_{011} V_{020} - B_{2h} V_{002} V_{020}\right]}{\overline{Y}\left[V_{020} V_{002} - V_{011}^2\right]}
\end{aligned}
\right\} \tag{3.6}
$$

Putting optimum values of $m_1$ and $m_2$ from (3.6), we obtained min MSE of proposed estimator $t_p$.

## 4. Efficiency comparison

In this section, the conditions for which the proposed estimator $t_p$ is better than $\overline{y}_{st}$, $t_1$, $t_2$, $t_3$, $t_4$, $t_5$, $t_6$, and $t_7$.

The variance is given by

$$V(\bar{y}_{st}) = \bar{Y}^2 V_{200} \tag{4.1}$$

To compare the efficiency of the proposed estimator with the existing estimator, from (4.1) and (2.3), (2.4), (2.9), (2.10), (2.11), (2.12) and (2.14), we have

$$V(\bar{y}_{st}) - MSE(t_p) = \bar{Y}^2 P_1 + P_2 - \bar{Y} P_3 \geq 0 \tag{4.2}$$

$$MSE(t_1) - MSE(t_p) = \bar{Y}^2 \left[ V_{020} - 2V_{110} \right] - \bar{Y}^2 P_1 - P_2 + \bar{Y} P_3 \geq 0 \tag{4.3}$$

$$MSE(t_2) - MSE(t_p) = \bar{Y}^2 \left[ \frac{V_{020}}{4} - V_{110} \right] - \bar{Y}^2 P_1 - P_2 + \bar{Y} P_3 \geq 0 \tag{4.4}$$

$$MSE(t_3) - MSE(t_p) = \bar{Y}^2 \left[ \frac{V_{020}}{4} + \frac{V_{002}}{4} - V_{110} - V_{101} + \frac{V_{011}}{2} \right] - \bar{Y}^2 P_1 - P_2 + \bar{Y} P_3 \geq 0 \tag{4.5}$$

$$MSE(t_4) - MSE(t_p) = \bar{Y}^2 \left[ \frac{V_{020}}{4} + \frac{V_{002}}{4} + V_{110} + V_{101} + \frac{V_{011}}{2} \right] - \bar{Y}^2 P_1 - P_2 + \bar{Y} P_3 \geq 0 \tag{4.6}$$

$$MSE(t_5) - MSE(t_p) = \bar{Y}^2 \left[ \frac{V_{020}}{4} + \frac{V_{002}}{4} - V_{110} + V_{101} + \frac{V_{011}}{2} \right] - \bar{Y}^2 P_1 - P_2 + \bar{Y} P_3 \geq 0 \tag{4.7}$$

$$MSE(t_6) - MSE(t_p) = \bar{Y}^2 \left[ \frac{V_{020}}{4} + \frac{V_{002}}{4} + V_{110} - V_{101} + \frac{V_{011}}{2} \right] - \bar{Y}^2 P_1 - P_2 + \bar{Y} P_3 \geq 0 \tag{4.8}$$

Using (4.2) - (4.8), we conclude that the proposed estimator outperforms than the estimators considered in literature.

## 5. Empirical study

In this section, we use the data set in Koyuncu and Kadilar [4]. The population statistics are given in Table 3.2.1. In this data set, the study variable (Y) is the number of teachers, the first auxiliary variable (X) is the number of students, and the second auxiliary variable (Z) is the number of classes in both primary and secondary schools.

**Table 5.1**: Data Statistics of Population

| | | |
|---|---|---|
| $N_1 = 127$ | $N_2 = 117$ | $N_3 = 103$ |
| $N_4 = 170$ | $N_5 = 205$ | $N_6 = 201$ |
| $n_1 = 31$ | $n_2 = 21$ | $n_3 = 29$ |
| $n_4 = 38$ | $n_5 = 22$ | $n_6 = 39$ |
| $S_{y1} = 883.835$ | $S_{y2} = 644$ | $S_{y3} = 1033.467$ |
| $S_{y4} = 810.585$ | $S_{y5} = 403.654$ | $S_{y6} = 711.723$ |
| $\overline{Y}_1 = 703.74$ | $\overline{Y}_2 = 413$ | $\overline{Y}_3 = 573.17$ |
| $\overline{Y}_4 = 424.66$ | $\overline{Y}_5 = 267.03$ | $\overline{Y}_6 = 393.84$ |
| $S_{x1} = 30486.751$ | $S_{x2} = 15180.760$ | $S_{x3} = 27549.697$ |
| $S_{x4} = 18218.931$ | $S_{x5} = 8997.776$ | $S_{x6} = 23094.141$ |
| $\overline{X}_1 = 20804.59$ | $\overline{X}_2 = 9211.79$ | $\overline{X}_3 = 14309.30$ |
| $\overline{X}_4 = 9478.85$ | $\overline{X}_5 = 5569.95$ | $\overline{X}_6 = 12997.59$ |
| $S_{xy1} = 25237153.52$ | $S_{xy2} = 9747942.85$ | $S_{xy3} = 28294397.04$ |
| $S_{xy1} = 14523885.53$ | $S_{xy1} = 3393591.75$ | $S_{xy6} = 15864573.97$ |
| $\rho_{xy1} = 0.936$ | $\rho_{xy2} = 0.996$ | $\rho_{xy3} = 0.994$ |
| $\rho_{xy4} = 0.983$ | $\rho_{xy5} = 0.989$ | $\rho_{xy6} = 0.965$ |
| $S_{z1} = 555.5816$ | $S_{z2} = 365.4576$ | $S_{z3} = 612.9509281$ |
| $S_{z4} = 458.0282$ | $S_{z5} = 260.8511$ | $S_{z6} = 397.0481$ |
| $\overline{Z}_1 = 498.28$ | $\overline{Z}_2 = 318.33$ | $\overline{Z}_3 = 431.36$ |
| $\overline{Z}_4 = 498.28$ | $\overline{Z}_5 = 227.20$ | $\overline{Z}_6 = 313.71$ |
| $S_{yz1} = 480688.2$ | $S_{yz2} = 230092.8$ | $S_{yz1} = 623019.3$ |
| $S_{yz1} = 364943.4$ | $S_{yz1} = 101539$ | $S_{yz1} = 277696.1$ |
| $S_{xz1} = 15914648$ | $S_{xz2} = 5379190$ | $S_{xz3} = 164900674.56$ |
| $S_{xz4} = 8041254$ | $S_{xz5} = 2144057$ | $S_{xz1} = 8857729$ |
| $\rho_{yz1} = 0.978914$ | $\rho_{yz2} = 0.9762$ | $\rho_{yz3} = 0.983511$ |
| $\rho_{yz4} = 0.982958$ | $\rho_{yz5} = 0.964342$ | $\rho_{yz6} = 0.982689$ |

We have computed the pre relative efficiency (PRE) of different estimators of $\overline{Y}_{st}$ with respect to $\overline{y}_{st}$ and complied in table 5.2:

**Table 5.2: Percent Relative Efficiencies (PRE) of estimator**

| S.No. | Estimators | PRE'S |
|:---:|:---:|:---:|
| 1 | $\overline{y}_{st}$ | 100 |
| 2 | $t_1$ | 1029.46 |
| 3 | $t_2$ | 370.17 |
| 4 | $t_3$ | 2045.43 |
| 5 | $t_4$ | 27.94 |
| 6 | $t_5$ | 126.41 |
| 7 | $t_6$ | 77.21 |
| 8 | $t_7$ | 2360.54 |
| 9 | $t_p$ | 4656.35 |

## 6. Conclusion

In this paper, we proposed a new estimator for estimating unknown population mean of study variable using information on two auxiliary variables. Expressions for bias and MSE of the estimator are derived up to first degree of approximation. The proposed estimator is compared with usual mean estimator and other considered estimators. A numerical study is carried out to support the theoretical results. In the table 5.2, the proposed estimator performs better than the usual sample mean and other considered estimators.

# Some Ratio and Product Estimators Using Known Value of Population Parameters

**Rajesh Singh, Sachin Malik\* and Florentin Smarandache\*\***

\*Department of Mathematics, SRM University

Delhi NCR, Sonepat- 131029, India

Department of Statistics, Banaras Hindu University

Varanasi-221005, India

\*\*Department of Mathematics, University of New Mexico

Gallup, NM 87301, USA

\*Corresponding Author, sachinkurava999@gmail.com

**Abstract**

In the present article, we proposed a family of estimators for estimating population means using known value of some population parameters. Khoshnevisan et al. [1] proposed a general family of estimators for estimating population means using known value of some population parameter(s) which after some substitutions led to some ratio and product estimators initially proposed by Sisodia and Dwivedi [2], Singh and Tailor [3], Pandey and Dubey [4], Adewara et al. [5], yadav and Kadilar [6]. The present family of estimators provides us significant improvement over previous families in theory. An empirical study is carried out to judge the merit of the proposed estimator.

**Keywords**: Ratio Estimator, Product Estimator, Population Parameter, Efficiency, Mean Square Error.

## 1.    Introduction

The problem of estimating the population mean in the presence of an auxiliary variable has been widely discussed in finite population sampling literature. Ratio, product and difference methods of estimation are good examples in this context. Ratio method of estimation is quite effective when there is high positive correlation between study and auxiliary variables. On the other hand, if correlation is negative (high), the product method of estimation can be employed efficiently.

In recent years, a number of research papers on ratio-type, exponential ratio-type and regression-type estimators have appeared, based on different types of transformations. Some important contributions in this area are due to Singh and Tailor [3], Shabbir and Gupta [7,8], Kadilar and Cingi [9,10], Khosnevisan et. al.(2007).

Khoshnevisan et al. [1] defined their family of estimators as

$$t = \bar{y}[\frac{a\bar{X} + b}{\alpha(a\bar{x} + b) + (1-\alpha)(a\bar{X} + b)}]^g$$

where $a(\neq 0)$, b are either real numbers or the functions of the known parameters of the auxiliary variable x such as standard deviation ($\sigma_x$), Coefficient of Variation ($C_x$), Skewness ($\beta_1(x)$), Kurtosis ($\beta_2(x)$) and Correlation Coefficient ($\rho$).

(i). When $\alpha$=0, a=0=b, g=0, we have the mean per unit estimator, $t_0 = \bar{y}$ with

$$MSE(t_0) = (\frac{N-n}{Nn})\bar{Y}^2 C_y^2 \tag{1.1}$$

(ii). When $\alpha$=1, a=1, b=0, g=1, we have the usual ratio estimator, $t_1 = \bar{y}(\frac{\bar{X}}{\bar{x}})$ with

$$MSE(t_1) = (\frac{N-n}{Nn})\bar{Y}^2(C_y^2 + C_x^2 - 2\rho\rho_x C_y) \tag{1.2}$$

(iii). When $\alpha$=1, a=1, b=0, g=-1, we have the usual product estimator, $t_2 = \bar{y}(\frac{\bar{x}}{\bar{X}})$ with

$$MSE(t_2) = (\frac{N-n}{Nn})\bar{Y}^2(C_y^2 + C_x^2 + 2\rho\rho_x C_y) \tag{1.3}$$

(iv). When $\alpha$=1, a=1, b=$C_x$, g=1, we have Sisodia and Dwivedi [2] ratio estimator, $t_3 = \bar{y}(\frac{\bar{X} + C_x}{\bar{x} + C_x})$ with

$$MSE(t_3) = (\frac{N-n}{Nn})\bar{Y}^2(C_y^2 + (\frac{\bar{X}}{\bar{X} + C_x})^2 C_x^2 - 2(\frac{\bar{X}}{\bar{X} + C_x})\rho\rho_x C_y) \tag{1.4}$$

(v). When $\alpha$=1, a=1, b=$C_x$, g=-1

we have Pandey and Dubey [4] product estimator, $t_4 = \bar{y}(\dfrac{\bar{x} + C_x}{\bar{X} + C_x})$ with

$$MSE(t_4) = (\frac{N-n}{Nn})\bar{Y}^2(C_y^{\,2} + (\frac{\bar{X}}{\bar{X} + C_x})^2 C_x^{\,2} + 2(\frac{\bar{X}}{\bar{X} + C_x})\rho\rho_x C_y) \qquad (1.5)$$

(vi). When $\alpha=1$, $a=1$, $b=\rho$, $g=1$, we have Singh and Taylor [3] ratio estimator, $t_5 = \bar{y}(\dfrac{\bar{X} + \rho}{\bar{x} + \rho})$ with

$$MSE(t_5) = (\frac{N-n}{Nn})\bar{Y}^2(C_y^{\,2} + (\frac{\bar{X}}{\bar{X} + \rho})^2 C_x^{\,2} - 2(\frac{\bar{X}}{\bar{X} + \rho})\rho\rho_x C_y) \qquad (1.6)$$

(vii). When $\alpha=1$, $a=1$, $b=\rho$, $g=-1$, we have Singh and Taylor [3] product estimator, $t_6 = \bar{y}(\dfrac{\bar{x} + \rho}{\bar{X} + \rho})$ with

$$MSE(t_6) = (\frac{N-n}{Nn})\bar{Y}^2(C_y^{\,2} + (\frac{\bar{X}}{\bar{X} + \rho})^2 C_x^{\,2} + 2(\frac{\bar{X}}{\bar{X} + \rho})\rho\rho_x C_y) \qquad (1.7)$$

There are other ratio and product estimators from these families that are not inferred here but this paper will be limited to those ones that made use of Coefficient of Variation ($C_x$) and Correlation Coefficient ($\rho$) since the conclusion obtained here can also be inferred on all others that made use of other population parameters such as the standard deviation ($\sigma_x$), Skewness ($\beta_1(x)$) and Kurtosis ($\beta_2(x)$) in the same family.

## 2. On the Modified Ratio and Product Estimators.

Adopting Adewara (2006), Adewara et al. (2012) proposed the following estimators as

$$t^*_{\,1} = \bar{y}^*(\frac{\bar{X}}{\bar{x}^*}), \qquad (2.1)$$

$$t^*_{\,2} = \bar{y}^*(\frac{\bar{x}^*}{\bar{X}}), \qquad (2.2)$$

$$t^*_{\,3} = \bar{y}^*(\frac{\bar{X} + C_x}{\bar{x}^* + C_x}), \qquad (2.3)$$

$$t^{*}_{4} = \overline{y}^{*}(\frac{\overline{x}^{*} + C_{x}}{\overline{X} + C_{x}}),$$ (2.4)

$$t^{*}_{5} = \overline{y}^{*}(\frac{\overline{X} + \rho}{\overline{x}^{*} + \rho}) \text{ and}$$ (2.5)

$$t^{*}_{6} = \overline{y}^{*}(\frac{\overline{x}^{*} + \rho}{\overline{X} + \rho}),$$ (2.6)

Where $\overline{x}^{*}$ and $\overline{y}^{*}$ are the sample means of the auxiliary variables and variable of interest yet to be drawn with the relationships (i) $\overline{X} = f\overline{x} + (1-f)\overline{x}^{*}$ and (ii) $\overline{Y} = f\overline{y} + (1-f)\overline{y}^{*}$. Srivenkataramana and Srinath [12].

The Mean Square Errors of these estimators $t^{*}_{i}$, i = 1,2, …, 6 are as follows:

(i). $MSE(t^{*}_{1}) = (\frac{n}{N-n})^{2} MSE(t_{1})$ (2.7)

(ii). $MSE(t^{*}_{2}) = (\frac{n}{N-n})^{2} MSE(t_{2})$ (2.8)

(iii). $MSE(t^{*}_{3}) = (\frac{n}{N-n})^{2} MSE(t_{3})$ (2.9)

(iv). $MSE(t^{*}_{4}) = (\frac{n}{N-n})^{2} MSE(t_{4})$ (2.10)

(v). $MSE(t^{*}_{5}) = (\frac{n}{N-n})^{2} MSE(t_{5})$ (2.11)

(vi). $MSE(t^{*}_{6}) = (\frac{n}{N-n})^{2} MSE(t_{6})$ (2.12)

Following Adewara et al [5], Yadav and Kadilar [6] proposed some improved ratio and product estimators for estimating the population mean of the study variable as follows

$$\eta^{*}_{1} = k\overline{y}^{*}(\frac{\overline{X}}{\overline{x}^{*}}),$$ (2.13)

$$\eta^{*}_{2} = k\,\overline{y}^{*}(\frac{\overline{x}^{*}}{\overline{X}}),$$ (2.14)

$$\eta^*_3 = k \, \overline{y}^* \left( \frac{\overline{X} + C_x}{\overline{x}^* + C_x} \right), \tag{2.15}$$

$$\eta^*_4 = k \, \overline{y}^* \left( \frac{\overline{x}^* + C_x}{\overline{X} + C_x} \right), \tag{2.16}$$

$$\eta^*_5 = k \, \overline{y}^* \left( \frac{\overline{X} + \rho}{\overline{x}^* + \rho} \right) \tag{2.17}$$

$$\eta^*_6 = k \, \overline{y}^* \left( \frac{\overline{x}^* + \rho}{\overline{X} + \rho} \right), \tag{2.18}$$

The mean square error of these estimators $\eta^*_i$, i=1,2,…,6 are as follows

$$MSE(\eta^*_1) = \overline{Y}^2 \left[ h^2 \left( k_1^2 \lambda C_y^2 + \left\{ 3k_1^2 - 2k_1 \right\} \lambda C_x^2 - 2 \left\{ 2k_1^2 - k_1 \right\} \lambda C_{yx} \right) + \left\{ k_1 - 1 \right\}^2 \right] \tag{2.19}$$

$$MSE(\eta^*_2) = \overline{Y}^2 \left[ h^2 \left( k_2^2 \lambda C_y^2 + k_2^2 \lambda C_x^2 + 2 \left\{ 2k_1^2 - k_1 \right\} \lambda C_{yx} \right) + \left\{ k_2 - 1 \right\}^2 \right] \tag{2.20}$$

$$MSE(\eta^*_3) = \overline{Y}^2 \left[ h^2 \left( k_3^2 \lambda C_y^2 + \left\{ 3k_3^2 - 2k_3 \right\} v_1^2 \lambda C_x^2 - 2v_1 \left\{ 2k_3^2 - k_3 \right\} \lambda C_{yx} \right) + \left\{ k_3 - 1 \right\}^2 \right] \tag{2.21}$$

$$MSE(\eta^*_4) = \overline{Y}^2 \left[ h^2 \left( k_4^2 \lambda C_y^2 + k_4^2 v_1^2 \lambda C_x^2 + 2v_1 \left\{ 2k_4^2 - k_4 \right\} \lambda C_{yx} \right) + \left\{ k_4 - 1 \right\}^2 \right] \tag{2.22}$$

$$MSE(\eta^*_3) = \overline{Y}^2 \left[ h^2 \left( k_3^2 \lambda C_y^2 + \left\{ 3k_3^2 - 2k_3 \right\} v_2^2 \lambda C_x^2 - 2v_2 \left\{ 2k_5^2 - k_5 \right\} \lambda C_{yx} \right) + \left\{ k_5 - 1 \right\}^2 \right] \tag{2.23}$$

$$MSE(\eta^*_4) = \overline{Y}^2 \left[ h^2 \left( k_6^2 \lambda C_y^2 + k_4^2 v_2^2 \lambda C_x^2 + 2v_2 \left\{ 2k_6^2 - k_6 \right\} \lambda C_{yx} \right) + \left\{ k_6 - 1 \right\}^2 \right] \tag{2.24}$$

Where,

$$\lambda = \frac{N\text{-}n}{Nn}, h = \frac{n}{N-n}, C_y^2 = \frac{S_y^2}{\overline{Y}^2}, C_x^2 = \frac{S_x^2}{\overline{X}^2}, C_{yx} = \frac{S_{yx}}{\overline{Y}\,\overline{X}}, v_1 = \frac{\overline{X}}{\overline{X} + C_x}, v_1 = \frac{\overline{X}}{\overline{X} + C_x} \text{ and } \rho = \frac{S_{yx}}{S_y S_x}$$

And $k_1 = \dfrac{h^2 \left[ \lambda C_x^2 - \lambda C_{yx} \right] + 1}{h^2 \left[ 3C_x^2 \lambda - 4C_{yx}\lambda + \lambda C_y^2 \right] + 1}, k_2 = \dfrac{h^2 \lambda C_{yx} + 1}{h^2 \left[ C_x^2 \lambda + 4C_{yx}\lambda + \lambda C_y^2 \right] + 1},$

$$k_3 = \frac{h^2 \left[ \lambda v_1^2 C_x^2 - v_1 \lambda C_{yx} \right] + 1}{h^2 \left[ 3v_1^2 C_x^2 \lambda - 4v_1 C_{yx}\lambda + \lambda C_y^2 \right] + 1}, k_4 = \frac{h^2 \left[ \lambda v_1 \lambda C_{yx} \right] + 1}{h^2 \left[ 3v_1^2 C_x^2 \lambda + 4v_1 C_{yx}\lambda + \lambda C_y^2 \right] + 1}$$

$$k_5 = \frac{h^2 \left[ \lambda v_2^2 C_x^2 - v_2 \lambda C_{yx} \right] + 1}{h^2 \left[ 3v_2^2 C_x^2 \lambda - 4v_2 C_{yx}\lambda + \lambda C_y^2 \right] + 1}, \text{and } k_6 = \frac{h^2 \left[ \lambda v_2 \lambda C_{yx} \right] + 1}{h^2 \left[ 3v_2^2 C_x^2 \lambda + 4v_2 C_{yx}\lambda + \lambda C_y^2 \right] + 1}$$

## 3. The Proposed family of estimators

Following Malik Singh [14], we define the following class of estimators for population mean $\overline{Y}$ as

$$t_M = \left\{ m_1 \overline{y}^* + m_2 \left[ \overline{X} - \overline{x}^* \right] \right\} \left( \frac{\psi \overline{X} + \delta}{\psi \overline{x}^* + \delta} \right)^\alpha \exp \left[ \frac{(\omega \overline{X} + \mu) - (\omega \overline{x}^* + \mu)}{(\omega \overline{X} + \mu) + (\omega \overline{x}^* + \mu)} \right]^\beta \tag{3.1}$$

Where $m_1$ and $m_2$ are suitably chosen constants. $\psi, \delta, \omega,$ and $\mu$ are either real numbers or function of known parameters of the auxiliary variable. The scalar $\alpha$ and $\beta$ takes values +1 and -1 for ratio and product type estimators respectively.

To obtain the MSE , let us define

$$\overline{y} = \overline{Y}(1 + e_0) \ , \overline{x} = \overline{X}(1 + e_1)$$

such that $\quad E(e_i) = 0 \ ,$ i=0,1 and

$$E(e_0^2) = \lambda C_y^2 , \ \ E(e_1^2) = \lambda C_x^2 , \ \ \ E(e_0 e_1) = \lambda \rho C_y C_x$$

expressing equation (3.1) in terms of e's and retaining only terms up to second degree of e's, we have

$$t_M = \left[ m_1 \overline{Y}(1 - h e_0) + m_2 \overline{X} h e_1 \right] \left\{ \frac{\psi \overline{X} + \delta}{\psi \overline{X}(1 - h e_1) + \delta} \right\}^\alpha \exp \left\{ \frac{\omega \overline{X} h e_1}{2 \omega \overline{X} + 2\mu - \omega \overline{X} h e_1} \right\}^\beta$$

$$= \left[ m_1 \overline{Y}(1 - h e_0) + m_2 \overline{X} h e_1 \right] \left\{ 1 - R_1 h e_1 \right\}^{-\alpha} \exp \left\{ \beta R_2 h e_1 \left( 1 - \frac{R_2 h e_1}{2} + \frac{R_2^2 h^2 e_1^2}{4} \right) \right\}$$

$$= m_1 \overline{Y} \left[ \begin{array}{c} 1 + \alpha h e_1 + \dfrac{\alpha(\alpha+1) h^2 e_1^2}{2} + \dfrac{\beta h e_1}{2} + \dfrac{\alpha \beta h^2 e_1^2}{2} + \dfrac{\beta^2 h^2 e_1^2}{8} + \dfrac{\beta h^2 e_1^2}{4} - h e_0 \\ - \alpha h^2 e_0 e_1 - \dfrac{\beta h^2 e_0 e_1}{2} \end{array} \right]$$

$$+ m_2 \overline{X} \left[ h e_1 + \alpha h^2 e_1^2 + \frac{\beta h^2 e_1^2}{2} \right] \tag{3.2}$$

where,     $R_1 = \dfrac{\psi\overline{X}}{\psi\overline{X}+\delta}, R_2 = \dfrac{\omega\overline{X}}{\omega\overline{X}+\mu}$

Subtracting $\overline{Y}$ from both the sides of (3.2), we have

$$\left(t_M - \overline{Y}\right) = m_1\overline{Y}\left[1 - he_0 + L_1e_1 + L_2e_1^2 - L_3e_0e_1\right] + m_2\overline{X}\left[he_1 + L_4e_1^2\right] - \overline{Y} \tag{3.3}$$

where,

$$L_1 = \alpha R_1 h + \frac{\beta h R_2}{2}$$

$$L_2 = \frac{\alpha(\alpha+1)h^2R_1^2}{2} + \frac{\alpha\beta h^2R_1R_2}{2} + \frac{\beta^2 h^2R_2^2}{8} + \frac{\beta h^2R_2^2}{4}$$

$$L_3 = \alpha h^2R_1 + \frac{\beta h^2R_2}{2}$$

$$L_4 = \alpha h^2R_1 + \frac{\beta h^2R_2}{2}$$

Squaring both sides of (3.3) and neglecting terms of e's having power greater than two, we have

$$MSE\left(t_M\right) = \overline{Y}^2\left[1 + m_1^2T_1 + m_2^2T_2 + 2m_1m_2T_3 - 2m_1T_4 - 2m_2T_5\right] \tag{3.4}$$

where,

$$T_1 = \overline{Y}^2\left[1 + \lambda h^2C_y^2 + L_1^2\lambda C_x^2 - 2hL_1\lambda\rho C_yC_x + 2L_2\lambda C_x^2 - 2L_3\lambda\rho C_yC_x\right]$$

$$T_2 = h^2\lambda\overline{X}^2C_x^2$$

$$T_3 = \overline{Y}\overline{X}\left[L_4\lambda C_x^2 + L_1\lambda hC_x^2 - h^2\lambda\rho C_yC_x\right]$$

$$T_4 = \overline{Y}^2\left[1 + L_2\lambda C_x^2 - L_3\lambda\rho C_yC_x\right]$$

$$T_5 = \overline{Y}\overline{X}L_4\lambda C_x^2$$

minimization of (3.4) with respect to $m_1$ and $m_2$ yields optimum values as

$$m_1 = \frac{\left(T_2T_4 - T_3T_5\right)}{T_1T_2 - T_3^2}, \quad m_2 = \frac{\left(T_1T_5 - T_3T_4\right)}{T_1T_2 - T_3^2}$$

## 4. Empirical Study:

## Population I: Kadilar and Cingi [9]

$N = 106$, $n = 20$,   $\rho = 0.86, C_y = 5.22, C_x = 2.1, \overline{Y} = 2212.59$ and $\overline{X} = 27421.70$

**Population II: Maddala [13]**

$N = 16$, $n = 4$, $\rho = -0.6823$, $C_y = 0.2278$, $C_x = 0.0986$, $\overline{Y} = 7.6375$ and $\overline{X} = 75.4313$

**4. Results:**

**Table 4.1:** Showing the estimates obtained for both the Khoshnevisan et al. [1] estimators and Adewara et al. [5] estimators

| Estimator | Population I ($\rho > 0$) | Population II ($\rho < 0$) |
|---|---|---|
| $t_0$ | 5411349 | 0.5676 |
| $t_1$ | 2542740 | - |
| $t_2$ | - | 0.3387 |
| $t_3$ | 2542893 | - |
| $t_4$ | - | 0.3388 |
| $t_5$ | 2542803 | - |
| $t_6$ | - | 0.3376 |
| $t_1^*$ | 137519.8 | - |
| $t_2^*$ | - | 0.03763 |
| $t_3^*$ | 137528 | - |
| $t_4^*$ | - | 0.03765 |
| $t_5^*$ | 137523.1 | - |
| $t_6^*$ | - | 0.03751 |

**Table 4.2:** Showing the estimates obtained for Yadav and Kadilar [6] estimators

| Estimator | Population I ($\rho > 0$) | Population II ($\rho < 0$) |
|---|---|---|
| $\eta^*_1$ | 136145.37 | - |
| $\eta^*_2$ | - | 0.03762 |
| $\eta^*_3$ | 136138.05 | - |
| $\eta^*_4$ | - | 0.03764 |
| $\eta^*_5$ | 136107.94 | - |
| $\eta^*_6$ | - | 0.03750 |

**Table 4.3**: MSE of suggested estimators with different values of constants

| | | | | | | | | | MSE | |
|---|---|---|---|---|---|---|---|---|---|---|
| $m_1$ | $m_2$ | $\alpha$ | $\beta$ | $\psi$ | $\delta$ | $\omega$ | $\mu$ | estimator | PopI | PopII |
| 1 | 0 | 1 | 0 | 1 | 0 | - | - | $t^*_1$ | 137519.8 | - |
| 1 | 0 | -1 | 0 | 1 | 0 | - | - | $t^*_2$ | - | 0.03763 |
| 1 | 0 | 1 | 0 | 1 | $C_x$ | - | - | $t^*_3$ | 137528 | - |
| 1 | 0 | -1 | 0 | 1 | $C_x$ | - | - | $t^*_4$ | - | 0.03765 |
| 1 | 0 | 1 | 0 | 1 | $\rho$ | - | - | $t^*_5$ | 137523.1 | - |
| 1 | 0 | -1 | 0 | 1 | $\rho$ | - | - | $t^*_6$ | - | 0.03751 |
| $m_1$ | 0 | 1 | 0 | 1 | 0 | - | - | $\eta^*_1$ | 136145.37 | - |
| $m_1$ | 0 | -1 | 0 | 1 | 0 | - | - | $\eta^*_2$ | - | 0.03762 |
| $m_1$ | 0 | 1 | 0 | 1 | $C_x$ | - | - | $\eta^*_3$ | 136138.05 | - |
| $m_1$ | 0 | -1 | 0 | 1 | $C_x$ | - | - | $\eta^*_4$ | - | 0.03764 |
| $m_1$ | 0 | 1 | 0 | 1 | $\rho$ | - | - | $\eta^*_5$ | 136107.94 | - |
| $m_1$ | 0 | -1 | 0 | 1 | $\rho$ | - | - | $\eta^*_6$ | - | 0.03750 |
| $m_1$ | $m_2$ | 1 | 1 | 1 | 1 | 1 | 1 | $t_M$ | 75502.23 | - |
| $m_1$ | $m_2$ | -1 | -1 | 1 | 1 | 1 | 1 | $t_M$ | - | 0.03370 |

Since conventionally, for ratio estimators to hold, $\rho > 0$ and also for product estimators to hold, $\rho < 0$. Therefore two data sets are used in this paper, one to determine the efficiency of the modified ratio estimators and the other to determine that of the product estimators as stated below.

## 5. Conclusion

In this paper, we have proposed a new family of estimator for estimating unknown population mean of study variable using auxiliary variable. Expressions for the MSE of the estimator are derived up to first order of approximation. The proposed family of estimator is compared with the several existing estimators in literature. From table 4.3, we observe that the new family of estimators performs better than the other estimators considered in this paper for both of the data sets.

# An Unbiased Estimator for Estimating Population Mean in Simple Random Sampling Using Auxiliary Attribute

**Rajesh Singh and Sachin Malik***

[*]Department of Mathematics, SRM University

Delhi NCR, Sonepat- 131029, India

Department of Statistics, Banaras Hindu University

Varanasi-221005, India

*Corresponding Author, sachinkurava999@gmail.com

**Abstract**

In this paper we have proposed an almost unbiased estimator using known value of some population parameter(s) with known population proportion of an auxiliary variable. A class of estimators is defined which includes Naik and Gupta [1], Singh and Solanki [2] and Sahai and Ray [3] estimators. Under simple random sampling without replacement (SRSWOR) scheme the expressions for bias and mean square error (MSE) are derived. Numerical illustrations are given in support of the present study.

**Key words**: Auxiliary information, bias, mean square error, unbiased estimator.

**Introduction**

It is well known that the precision of the estimates of the population mean or total of the study variable y can be considering improved by the use of known information on an auxiliary variable x which is highly correlated with the study variable y. Out of many methods ratio, product and regression methods of estimation are good illustrations in this context. Using known values of certain populations parameters several authors have proposed improved estimators including Singh and Tailor [4], Kadilar and Cingi [5], Gupta and Shabbir [6,7], Khoshnevisan et al. [8], Singh et al. [9], Singh et al. [10], Koyuncu and Kadilar [11], Diana et al. [12], Upadhyaya et al. [13] and Singh and Solanki [2].

In many practical situations, instead of existence of auxiliary variables there exit some auxiliary attributes $\phi$ (say), which are highly correlated with the study variable y, such as

i.  Amount of milk produced (y) and a particular breed of cow ($\phi$).

ii.  Sex ($\phi$) and height of persons (y) and

iii.  Amount of yield of wheat crop and a particular variety of wheat ($\phi$) etc. (see Jhajj at al. [15]).

Many more situations can be encountered in practice where the information of the population mean $\overline{Y}$ of the study variable y in the presence of auxiliary attributes assumes importance. For these reasons various authors such as Naik and Gupta [1], Jhajj et al. [14], Abd- Elfattah et al. [15], Grover and Kaur [16], Malik and Singh [17] and Singh and Solanki [2] have paid their attention towards the improved estimation of population mean $\overline{Y}$ of the study variable y taking into consideration the point bi-serial correlation between a variable and an attribute.

Let $A = \sum_{i=1}^{N} \phi_i$ and $a = \sum_{i=1}^{n} \varphi_i$ denote the total number of units in the population and sample possessing attribute $\phi$ respectively, $p = \dfrac{A}{N}$ and $p = \dfrac{a}{n}$ denote the proportion of units in the population and sample, respectively, possessing attribute $\phi$.

Define,

$$e_y = \frac{\left(\overline{y} - \overline{Y}\right)}{\overline{Y}} \quad e_\phi = \frac{\left(p - P\right)}{P},$$

$$E\left(e_i\right) = 0, \left(i = y, \phi\right)$$

$$E\left(e_y^2\right) = fC_y^2, \quad E\left(e_\phi^2\right) = fC_p^2, \quad E\left(e_y e_\phi\right) = f\rho_{pb}C_y C_{p.}$$

Where

$$f_1 = \left(\frac{1}{n} - \frac{1}{N}\right), \quad C_y^2 = \frac{S_y^2}{\overline{Y}^2}, \quad C_p^2 = \frac{S_p^2}{P^2},$$

and $\rho_{pb} = \dfrac{S_{y\phi}}{S_y S_\phi}$ is the point bi-serial correlation coefficient.

Here,

$$S_y^2 = \frac{1}{N-1}\sum_{i=1}^{N}\left(y_i - \overline{Y}\right)^2, \, S_\phi^2 = \frac{1}{N-1}\sum_{i=1}^{N}\left(\phi_i - P\right)^2 \text{ and } S_{y\phi} = \frac{1}{N-1}\left(\sum_{i=1}^{N} y_i \phi_i - NP\overline{Y}\right)$$

In order to have an estimate of the study variable y, assuming the knowledge of the population proportion P, Naik and Gupta [1] proposed following estimate

$$t_{NGR} = \bar{y}\left(\frac{P}{p}\right)$$

(1.1)

$$t_{NGP} = \bar{y}\left(\frac{p}{P}\right)$$

(1.2)

following Naik and Gupta [1] , we propose the following estimator

$$t_1 = \bar{y}\left(\frac{K_1 P + K_2 K_3}{K_1 p + K_2 K_3}\right)^{\alpha}$$

(1.3)

The Bias and MSE expression's of the estimator $t_1$ up to the first order of approximation are, respectively, given by

$$B(t_1) = \bar{Y} f_1 C_p^2 \left[\frac{\alpha(\alpha+1)V_1^2}{2} - \alpha V_1 K_p\right]$$

(1.4)

$$MSE(t_1) = \bar{Y}^2 f_1 \left[C_y^2 + C_p^2\left(\alpha^2 V_1^2 - 2\alpha V_1 K_{p_1}\right)\right]$$

(1.5)

Also following Singh and Solanki [2], we propose the following estimator

$$t_2 = \bar{y}\left\{2 - \left(\frac{p}{P}\right)^{\beta} \exp\left[\lambda\left(\frac{(K_4 P + K_5) - (K_4 p + K_5)}{(K_4 P + K_5) + (K_4 p + K_5)}\right)\right]\right\}$$

(1.6)

The Bias and MSE expression's of the estimator $t_2$ up to the first order of approximation are, respectively, given by

$$B(t_2) = \bar{Y} f_1 C_p^2 \left[\frac{\lambda V_2 \beta}{2} - \frac{\beta(\beta-1)}{2} - \frac{\lambda(\lambda+2)V_2^2}{8} - \beta K_p + \frac{\lambda V_2 K_p}{2}\right]$$

(1.7)

$$MSE(t_2) = \bar{Y}^2 f_1 \left[C_y^2 + C_p^2\left(\beta^2 + \frac{\lambda^2 V_2^2}{4} - \beta\lambda V_2\right) - 2K_p C_p^2\left(\beta - \frac{\lambda V_2}{2}\right)\right]$$

(1.8)

$\alpha$, $\lambda$ and $\beta$ are suitable chosen constants. Also $K_1, K_3, K_4, K_5$ are either real numbers or function of known parameters of the auxiliary attributes $\phi$ such as $C_p$, $\beta_2(\phi)$, $\rho_{pb}$ and $K_p$. $K_2$ is an integer which takes values +1 and -1 for designing the estimators and

$$\left.\begin{aligned} V_1 &= \frac{K_1 P}{K_1 P + K_2 K_3} \\ V_2 &= \frac{K_4 P}{K_4 P + K_5} \end{aligned}\right\}$$

We see that the estimators $t_1$ and $t_2$ are biased estimators. In some applications bias is disadvantageous. Following these estimators we have proposed almost unbiased estimator of $\overline{Y}$.

## 2. Almost unbiased estimator

Suppose $t_0 = \overline{y}$, $t_1 = \overline{y}\left(\dfrac{K_1 P + K_2 K_3}{K_1 p + K_2 K_3}\right)^\alpha$, $t_2 = \overline{y}\left\{2 - \left(\dfrac{p}{P}\right)^\beta \exp\left[\lambda\left(\dfrac{(K_4 P + K_5) - (K_4 p + K_5)}{(K_4 P + K_5) + (K_4 p + K_5)}\right)\right]\right\}$

Such that $t_0$, $t_1$, $t_2 \in W$, where W denotes the set of all possible estimators for estimating the population mean $\overline{Y}$. By definition, the set W is a linear variety if

$$t_p = \sum_{i=0}^{3} w_i t_i \in W \tag{2.1}$$

Such that, $\displaystyle\sum_{i=0}^{3} w_i = 1$ and $w_i \in R$ \hfill (2.2)

where $w_i (i = 0,1,2,3)$ denotes the constants used for reducing the bias in the class of estimators, H denotes the set of those estimators that can be constructed from $t_i (i = 0,1,2,3)$ and R denotes the set of real numbers.

Expressing $t_p$ in terms of e's, we have

$$t_p = \overline{Y}\left[1 + e_y + w_1\left(\frac{\alpha(\alpha+1)V_1^2 e_\phi^2}{2} - \alpha V_1 e_\phi - \alpha V_1 e_y e_\phi\right)\right.$$

$$+ w_2 \left( -\beta e_\phi - \frac{\beta(\beta-1)e_\phi^2}{2} + \frac{\lambda V_2 e_\phi}{2} + \frac{\lambda V_2 \beta e_\phi^2}{2} - \frac{\lambda(\lambda+2)V_2^2 e_\phi^2}{8} - \beta e_y e_\phi + \frac{\lambda V_2 e_y e_\phi}{2} \right) \right] \qquad (2.3)$$

Subtracting $\overline{Y}$ from both sides of equation (2.3) and then taking expectation of both sides, we get the bias of the estimator $t_6$ up to the first order of approximation, as

$$B(t_p) = \overline{Y} f_1 w_1 C_p^2 \left( \frac{\alpha(\alpha+1)V_1^2}{2} - \alpha V_1 K_p \right) + \overline{Y} f_1 w_2 C_p^2 \left( \frac{\lambda V_2 \beta}{2} - \frac{\beta(\beta-1)}{2} - \frac{\lambda(\lambda+2)V_2^2}{8} \right.$$

$$\left. - \beta K_p + \frac{\lambda V_2 K_p}{2} \right) \qquad (2.4)$$

From (2.3), we have

$$(t_p - \overline{Y}) = \overline{Y} \left[ e_0 - w_1 \alpha V_1 e_\phi - w_2 \left( \beta e_\phi + \frac{\lambda V_2 e_\phi}{2} \right) \right]$$

$$(2.5)$$

Squaring both sides of (2.5) and then taking expectation, we get the MSE of the estimator $t_6$ up to the first order of approximation, as

$$MSE(t_p) = \overline{Y}^2 f_1 \left[ C_y^2 + C_p^2 \left( Q^2 - 2Q K_p \right) \right] \qquad (2.6)$$

Where

Which is minimum when

$$Q = K_p \qquad (2.7)$$

Where $Q = w_1 \alpha V_1 + w_2 \left( \beta - \frac{\lambda V_2}{2} \right) \qquad (2.8)$

Putting the value of $Q = K_p$ in (2.6) we have optimum value of estimator as $t_p$ (optimum).

Thus the minimum MSE of $t_p$ is given by

$$\min.MSE(t_p) = \overline{Y}^2 f_1 C_y^2 \left( 1 - \rho_{pb}^2 \right) \qquad (2.9)$$

Which is same as that of traditional linear regression estimator.

from (2.2) and (2.8), we have only two equations in three unknowns. It is not possible to find the unique values for $w_i$'s, $1=0,1,2$. In order to get unique values of $w_i$'s, we shall impose the linear restriction

$$\sum_{i=0}^{2} w_i B(t_i) = 0 \qquad (2.10)$$

where $B(t_i)$ denotes the bias in the i$^{th}$ estimator.

Equations (2.2), (2.8) and (2.10) can be written in the matrix form as

$$
\begin{bmatrix}
1 & 1 & 1 \\
0 & \alpha V_1 & \beta - \dfrac{\lambda V_2}{2} \\
0 & B(t_1) & B(t_2)
\end{bmatrix}
\begin{bmatrix}
w_0 \\
w_1 \\
w_2
\end{bmatrix}
=
\begin{bmatrix}
1 \\
k_p \\
0
\end{bmatrix}
\qquad (2.11)
$$

Using (2.11), we get the unique values of $w_i$'s, $1=0,1,2$ as

$$
\left.
\begin{aligned}
w_0 &= \frac{\alpha V_1[\alpha V_1 A_2 - A_1 X_1] - X_1 K_P A_1 - X_2 \alpha V_1[\alpha V_1 A_2 - A_1 X_1] - \alpha V_1 K_P A_1}{\alpha V_1[\alpha V_1 A_2 - A_1 X_1]} \\
w_1 &= \frac{X_1 K_P A_1}{\alpha V_1[\alpha V_1 A_2 - A_1 X_1]} + X_2 \\
w_2 &= \frac{K_P A_1}{[\alpha V_1 A_2 - A_1 X_1]}
\end{aligned}
\right\}
$$

where,

$$
\left.
\begin{aligned}
A_1 &= \frac{\alpha(\alpha+1)V_1^2}{2} - \alpha V_1 K_P \\
A_2 &= \frac{\lambda V_2 \beta}{2} - \frac{\beta(\beta-1)}{2} - \frac{\lambda(\lambda+2)V_2^2}{8} - \beta K_P + \frac{\lambda V_2 K_P}{2} \\
X_1 &= A_1\left[\beta - \frac{\lambda V_2}{2}\right] \\
X_2 &= \frac{K_p}{\alpha V_1}
\end{aligned}
\right\}
$$

Use of these $w_i$'s, $1=0,1,2$ remove the bias up to terms of order $o(n^{-1})$ at (2.1).

## 3. Empirical study

For empirical study we use the data sets earlier used by Sukhatme and Sukhatme [18], (p.256) (population 1) and Mukhopadhyaya [19] (p.44) (population 2) to verify the theoretical results.

**Data statistics:**

| Population | N | n | $\overline{Y}$ | P | $C_y$ | $C_p$ | $\rho_{pb}$ | $\beta_2(\phi)$ |
|------------|----|----|-------|--------|---------|---------|--------|---------|
| **Population 1** | 89 | 20 | 3.360 | 0.1236 | 0.60400 | 2.19012 | 0.766 | 6.2381 |
| **Population 2** | 25 | 10 | 9.44 | 0.400 | 0.17028 | 1.27478 | -0.387 | 4.3275 |

**Table 3.1** : Values of $w_i$'s,

| $w_i$'s, | **Population 1** | **Population 2** |
|----------|------------------|------------------|
| $w_0$ | -3.95624 | 1.124182 |
| $w_1$ | 5.356173 | 0.020794 |
| $w_2$ | 0.39993 | -0.14498 |

**Table 3.2: PRE of different estimators of $\overline{Y}$ with respect to $\overline{y}$**

| Choice of scalars | | | | | | | | | | | Estimator | PRE (POPII) | PRE (POPII) |
|------|------|------|------|------|------|------|------|------|------|------|-----------|-------------|-------------|
| $w_0$ | $w_1$ | $w_2$ | $K_1$ | $K_2$ | $K_3$ | $K_4$ | $K_5$ | $\alpha$ | $\beta$ | $\lambda$ | | | |
| 1 | 0 | 0 | | | | | | | | | $\overline{y}$ | 100 | 100 |
| 0 | 1 | 0 | 1 | 1 | 0 | | 1 | | | | $t_{NGR}$ | 11.63 | 1.59 |
| | | | 1 | 1 | 0 | | -1 | | | | $t_{NGP}$ | 5.075 | 1.94 |
| | | | | | | | | | | | | | |
| 0 | 0 | 1 | | | | | | | 1 | 0 | $t_{1(1,0)}$ | 12.88 | 1.59 |
| | | | | | | | | | -1 | 0 | $t_{1(-1,0)}$ | 5.43 | 1.95 |
| | | | | | 1 | 0 | | | 1 | 1 | $t_{2(1,1)}$ | 73.59 | 0.84 |
| | | | | | 1 | 0 | | | 1 | -1 | $t_{2(1,-1)}$ | 4.94 | 8.25 |
| | | | | | 1 | 0 | | | 0 | 1 | $t_{2(0,1)}$ | 14.95 | 8.25 |
| | | | | | 1 | 0 | | | 0 | -1 | $t_{2(0,-1)}$ | 73.48 | 5.58 |
| $w_0$ | $w_1$ | $w_2$ | 1 | 1 | 1 | 1 | 1 | 1 | 1 | 1 | $t_P$ **optimum** | **241.98** | **117.61** |

## 4. Proposed estimators in two phase sampling

In some practical situations when P is not known a priori, the technique of two-phase sampling is used. Let $p'$ denote the proportion of units possessing attribute $\phi$ in the first phase sample of size $n'$; $p$ denote the proportion of units possessing attribute $\phi$ in the second phase sample of size $n' > n$ and $\overline{y}$ denote the mean of the study variable y in the second phase sample.

In two-phase sampling the estimator $t_p$ will take the following form

$$t_{pd} = \sum_{i=0}^{3} h_i t_{id} \in H \tag{4.1}$$

Such that, $\sum_{i=0}^{3} h_i = 1$ and $h_i \in R$ (4.2)

Where,

$$t_{0d} = \overline{y} \ , t_{1d} = \overline{y} \left( \frac{K_1 p' + K_2 K_3}{K_1 p + K_2 K_3} \right)^m , t_{2d} = \overline{y} \left\{ 2 - \left( \frac{p}{p'} \right)^n \exp \left[ \gamma \left( \frac{(K_4 p' + K_5) - (K_4 p + K_5)}{(K_4 p' + K_5) + (K_4 p + K_5)} \right) \right] \right\}$$

The Bias and MSE expression's of the estimator $t_{1d}$ and $t_{2d}$ up to the first order of approximation are, respectively, given by

$$B(t_{1d}) = \overline{Y}\left[\frac{m(m-1)R_1^2 f_2 C_p^2}{2} + \frac{m(m+1)R_1^2 f_1 C_p^2}{2} - m^2 R_1^2 f_2 C_p^2 + mR_1 f_3 k_p C_p^2\right] \qquad (4.3)$$

$$MSE(t_1) = \overline{Y}^2\left[f_1 C_y^2 + m^2 R_1^2 f_3 C_p^2 - 2mR_1 k_p f_3 C_p^2\right] \qquad (4.4)$$

$$B(t_{2d}) = \overline{Y}\left[-\frac{n(n-1)f_1 C_p^2}{2} + \frac{n(n+1)f_2 C_p^2}{2} + nf_2 k_p C_p^2 + n^2 f_2 C_p^2 + f_3 \gamma R_2 k_p C_p^2 + f_3 \gamma R_2 n C_p^2\right] \qquad (4.5)$$

$$MSE(t_{2d}) = \overline{Y}^2\left[f_1 C_y^2 + L_1^2 f_3 C_p^2\right] \qquad (4.6)$$

Where,

$$\left.\begin{array}{l} R_1 = \dfrac{K_1 P}{K_1 P + K_2 K_3} \\[3mm] R_2 = \dfrac{K_4 P}{2[K_4 P + K_5]} \\[3mm] L_1 = n - \gamma A_2 \end{array}\right\} \qquad (4.7)$$

Expressing (4.1) in terms of e's, we have

$$t_p = \overline{Y}\left[1 + e_y + w_1\left(\frac{m(m+1)R_1^2 e_\varphi^2}{2} - mR_1 e_\varphi - mR_1 e_y e_\varphi + mR_1 e_y e_\varphi + \frac{m(m-1)R_1^2 e_\varphi'^2}{2} + mR_1 e_y e_\varphi'\right)\right.$$

$$\left. + w_2\left(-ne_\varphi - \frac{n(n-1)e_\varphi^2}{2} + ne'_\varphi + n^2 e_\varphi e'_\varphi - \frac{n(n+1)e_\varphi'^2}{2} - \gamma R_2\left(e'_\varphi - e_\varphi\right) + \gamma R_2\left(e_y e_\varphi - e_y e'_\varphi\right) - ne_y e_\varphi\right)\right]$$

Subtracting $\overline{Y}$ from both sides of equation (2.3) and then taking expectation of both sides, we get the bias of the estimator $t_6$ up to the first order of approximation, as

$$B(t_p) = \overline{Y}\left[B(t_{1d}) + B(t_{2d})\right] \qquad (4.8)$$

Also,

$$(t_p - \overline{Y}) = \overline{Y}\left[e_0 + w_1\left[mR_1e'_\varphi - mR_1e_\varphi\right] + w_2\left(-ne_\varphi + ne'_\varphi - \gamma R_2e'_\varphi + \gamma R_2e_\varphi\right)\right]$$

$$(4.9)$$

Squaring both sides of (4.9) and then taking expectation, we get the MSE of the estimator $t_p$ up to the first order of approximation, as

$$MSE(t_p) = \overline{Y}^2 f_1 C_y^2 + L_2^2 f_3 C_p^2 - 2L_2 f_3 k_p C_p^2$$

$$(4.10)$$

Where

Which is minimum when

$$L_2 = K_p$$

$$(4.11)$$

Where $L_2 = w_1mR_1 + w_2(n - \gamma R_2)$

$$(4.12)$$

Putting the value of $L_2 = K_p$ in (4.10), we have optimum value of estimator as $t_p$ (optimum).

Thus the minimum MSE of $t_p$ is given by

$$\min.MSE(t_p) = \overline{Y}^2 C_y^2\left(f_1 - f_3\rho_{pb}^2\right)$$

$$(4.13)$$

Which is same as that of traditional linear regression estimator.

from (4.2) and (4.12), we have only two equations in three unknowns. It is not possible to find the unique values for $h_i$'s, $1=0,1,2$. In order to get unique values of $h_i$'s, we shall impose the linear restriction

$$\sum_{i=0}^{2} h_i B(t_i) = 0$$

$$(4.14)$$

where $B(t_i)$ denotes the bias in the $i^{th}$ estimator.

Equations (4.2), (4.12) and (4.14) can be written in the matrix form as

$$\begin{bmatrix} 1 & 1 & 1 \\ 0 & mR_1 & n - \gamma R_2 \\ 0 & B(t_1) & B(t_2) \end{bmatrix}\begin{bmatrix} h_0 \\ h_1 \\ h_2 \end{bmatrix} = \begin{bmatrix} 1 \\ k_p \\ 0 \end{bmatrix}$$

$$(4.15)$$

Using (4.15), we get the unique values of $w_i$'s, $1=0,1,2$ as

$$\left.\begin{array}{l} h_0 = 1 - h_1 - h_2 \\[2mm] h_1 = \dfrac{k_p}{mR_1} - \dfrac{P_1 K_P (n - \gamma R_2)}{P_1 n - m R_1 P_2 - P_1 \gamma R_2} \\[4mm] h_2 = \dfrac{K_P P_1}{\left[ P_1 n - m R_1 P_2 - P_1 \gamma R_2 \right]} \end{array}\right\}$$

Where,

$$\left.\begin{array}{l} P_1 = \dfrac{m(m-1)R_1^2 f_2 C_p^2}{2} + \dfrac{m(m+1)R_1^2 f_1 C_p^2}{2} - m^2 R_1^2 f_2 C_p^2 + m R_1 f_3 k_p C_p^2 \\[4mm] P_2 = \left[ -\dfrac{n(n-1)f_1 C_p^2}{2} + \dfrac{n(n+1)f_2 C_p^2}{2} + n f_2 k_p C_p^2 + n^2 f_2 C_p^2 + f_3 \gamma R_2 k_p C_p^2 + f_3 \gamma R_2 n C_p^2 \right] \end{array}\right\}$$

Use of these $h_i$'s, $1=0,1,2$ remove the bias up to terms of order $o\left(n^{-1}\right)$ at (4.1).

## 5. Empirical Study

For empirical study we use the data sets earlier used by Sukhatme and Sukhatme [18] (, p.256) (population 1) and  Mukhopadhyaya  [19] ( p.44) (population 2) to verify the theoretical results.

### 5.1 Data statistics:

| Pop. | N | n | $\overline{Y}$ | P | p' | $C_y$ | $C_p$ | $\rho_{pb}$ | n' |
|------|-----|----|------|--------|---------|---------|--------|--------|-----|
| **Pop. 1** | 89 | 23 | 1322 | 0.1304 | 0.13336 | 0.69144 | 2.7005 | 0.408 | 45 |
| **Pop. 2** | 25 | 7 | 7.143 | 0.294 | 0.308 | 0.36442 | 1.3470 | -0.314 | 13 |

**Table 5.2:  PRE of different estimators of $\bar{Y}$  with respect to $\bar{y}$**

| Choice of scalars | | | | | | | | | | | | | |
|---|---|---|---|---|---|---|---|---|---|---|---|---|---|
| $w_0$ | $w_1$ | $w_2$ | $K_1$ | $K_2$ | $K_3$ | $K_4$ | $K_5$ | $M$ | $n$ | $\gamma$ | Estimator | PRE (POPII) | PRE (POPII) |
| 1 | 0 | 0 | | | | | | | | | $\bar{y}$ | 100 | 100 |
| 0 | 1 | 0 | 1 | 1 | 0 | | | 1 | | | $t_{NGR}$ | 11.13 | 8.85 |
| | | | 1 | 1 | 0 | | | -1 | | | $t_{NGP}$ | 7.48 | 12.15 |
| 0 | 0 | 1 | | | | | | 1 | | 0 | $t_{1d(1,0)}$ | 26.84 | 5.42 |
| | | | | | | | | -1 | | 0 | $t_{1d(-1,0)}$ | 23.75 | 5.87 |
| | | | | | | 1 | 0 | | 1 | 1 | $t_{2d(1,1)}$ | 82.55 | 1.23 |
| | | | | | | 1 | 0 | | 1 | -1 | $t_{2d(1,-1)}$ | 8.56 | 8.46 |
| | | | | | | 1 | 0 | | 0 | 1 | $t_{2d(0,1)}$ | 22.54 | 6.57 |
| | | | | | | 1 | 0 | | 0 | -1 | $t_{2d(0,-1)}$ | 82.56 | 7.45 |
| $w_0$ | $w_1$ | $w_2$ | 1 | 1 | 1 | 1 | 1 | 1 | 1 | 1 | $t_P$ **optimum** | 112.55 | 106.89 |

## Conclusion

In this paper, we have proposed an unbiased estimator $t_p$ and $t_{pd}$ using information on the auxiliary attribute(s) in case of single phase and double phase sampling respectively. Expressions for bias and MSE's of the proposed estimators are derived up to first degree of approximation.

From theoretical discussion and empirical study we conclude that the proposed estimators $t_p$ and $t_{pd}$ under optimum conditions perform better than other estimators considered in the article.

**Appendix A.**

**Some members of the proposed family of estimators -**

| Some members (ratio-type) of the class $t_1$ When $w_0 = 0$, $w_1 = 1$, $w_2 = 0$ $\alpha = 1$, | | | | | |
|---|---|---|---|---|---|
| $K_1$ | $K_3$ | Estimators ($K_2 = 1$) | Estimators ($K_2 = -1$) | PRE'S $K_2 = 1$ | PRE'S $K_2 = -1$ |
| 1 | $C_p$ | $t_{1a1} = \bar{y}\left[\dfrac{P + C_p}{p + C_p}\right]$ | $t_{1b1} = \bar{y}\left[\dfrac{P - C_p}{p - C_p}\right]$ | 134.99 | 72.50 |
| 1 | $\beta_2(\phi)$ | $t_{1a2} = \bar{y}\left[\dfrac{P + \beta_2(\phi)}{p + \beta_2(\phi)}\right]$ | $t_{1b2} = \bar{y}\left[\dfrac{P - \beta_2(\phi)}{p - \beta_2(\phi)}\right]$ | 111.62 | 89.34 |
| $\beta_2(\phi)$ | $C_p$ | $t_{1a3} = \bar{y}\left[\dfrac{P\beta_2(\phi) + C_p}{p\beta_2(\phi) + C_p}\right]$ | $t_{1b3} = \bar{y}\left[\dfrac{P\beta_2(\phi) - C_p}{p\beta_2(\phi) - C_p}\right]$ | 226.28 | 12.99 |
| $C_p$ | $\beta_2(\phi)$ | $t_{1a4} = \bar{y}\left[\dfrac{PC_p + \beta_2(\phi)}{pC_p + \beta_2(\phi)}\right]$ | $t_{1b4} = \bar{y}\left[\dfrac{PC_p - \beta_2(\phi)}{pC_p - \beta_2(\phi)}\right]$ | 126.66 | 77.93 |
| 1 | $\rho_{pb}$ | $t_{1a5} = \bar{y}\left[\dfrac{P + \rho_{pb}}{p + \rho_{pb}}\right]$ | $t_{1b5} = \bar{y}\left[\dfrac{P - \rho_{pb}}{p - \rho_{pb}}\right]$ | 207.46 | 39.13 |
| NP | $S_\phi$ | $t_{1a6} = \overline{y}\left[\dfrac{NP^2 + S_\phi}{NPp + S_\phi}\right]$ | $t_{1b6} = \overline{y}\left[\dfrac{NP^2 - S_\phi}{NPp - S_\phi}\right]$ | 18.14 | 6.86 |
| NP | f | $t_{1a7} = \overline{y}\left[\dfrac{NP^2 + f}{NPp + f}\right]$ | $t_{1b7} = \overline{y}\left[\dfrac{NP^2 - f}{NPp - f}\right]$ | 13.79 | 10.85 |
| $\beta_2(\phi)$ | $K_{pb}$ | $t_{1a8} = \overline{y}\left[\dfrac{\beta_2(\phi)P + K_{pb}}{\beta_2(\phi)p + K_{pb}}\right]$ | $t_{1b8} = \overline{y}\left[\dfrac{\beta_2(\phi)P - K_{pb}}{\beta_2(\phi)p - K_{pb}}\right]$ | 24.15 | 5.40 |
| NP | $K_{pb}$ | $t_{1a9} = \overline{y}\left[\dfrac{NP^2 + K_{pb}}{NPp + K_{pb}}\right]$ | $t_{1b9} = \overline{y}\left[\dfrac{NP^2 - K_{pb}}{NPp - K_{pb}}\right]$ | 18.62 | 7.78 |
| N | 1 | $t_{1a10} = \overline{y}\left[\dfrac{NP + 1}{Np + 1}\right]$ | $t_{1b10} = \overline{y}\left[\dfrac{NP - 1}{Np - 1}\right]$ | 15.93 | 9.26 |
| N | $C_p$ | $t_{1a11} = \overline{y}\left[\dfrac{NP + C_p}{Np + C_p}\right]$ | $t_{1b11} = \overline{y}\left[\dfrac{NP - C_p}{Np - C_p}\right]$ | 19.79 | 6.86 |
| N | $\rho_{pb}$ | $t_{1a12} = \overline{y}\left[\dfrac{NP + \rho_{pb}}{Np + \rho_{pb}}\right]$ | $t_{1b12} = \overline{y}\left[\dfrac{NP - \rho_{pb}}{Np - \rho_{pb}}\right]$ | 15.18 | 9.78 |
| N | $S_\phi$ | $t_{1a13} = \overline{y}\left[\dfrac{NP + S_\phi}{Np + S_\phi}\right]$ | $t_{1b13} = \overline{y}\left[\dfrac{NP - S_\phi}{Np - S_\phi}\right]$ | 12.34 | 10.96 |
| N | f | $t_{1a14} = \overline{y}\left[\dfrac{NP + f}{Np + f}\right]$ | $t_{1b14} = \overline{y}\left[\dfrac{NP - f}{Np - f}\right]$ | 12.99 | 11.54 |
| N | g=1-f | $t_{1a15} = \overline{y}\left[\dfrac{NP + g}{Np + g}\right]$ | $t_{1b15} = \overline{y}\left[\dfrac{NP - g}{Np - g}\right]$ | 15.81 | 9.34 |
| N | $K_{pb}$ | $t_{1a16} = \overline{y}\left[\dfrac{NP + K_{pb}}{Np + K_{pb}}\right]$ | $t_{1b16} = \overline{y}\left[\dfrac{NP - K_{pb}}{Np - K_{pb}}\right]$ | 13.52 | 11.10 |
| n | $\rho_{pb}$ | $t_{1a17} = \bar{y}\left[\dfrac{nP + \rho_{pb}}{np + \rho_{pb}}\right]$ | $t_{1b17} = \bar{y}\left[\dfrac{nP - \rho_{pb}}{np - \rho_{pb}}\right]$ | 25.13 | 4.86 |
| n | $S_\phi$ | $t_{1a18} = \bar{y}\left[\dfrac{nP + \rho_{pb}}{np + \rho_{pb}}\right]$ | $t_{1b18} = \bar{y}\left[\dfrac{nP - \rho_{pb}}{np - \rho_{pb}}\right]$ | 14.98 | 8.81 |
| n | f | $t_{1a19} = \bar{y}\left[\dfrac{nP + f}{np + f}\right]$ | $t_{1b19} = \bar{y}\left[\dfrac{nP - f}{np - f}\right]$ | 13.38 | 11.20 |
| n | g=1-f | $t_{1a20} = \bar{y}\left[\dfrac{nP + g}{np + g}\right]$ | $t_{1b20} = \bar{y}\left[\dfrac{nP - g}{np - g}\right]$ | 29.13 | 3.68 |
| n | $K_{pb}$ | $t_{1a21} = \bar{y}\left[\dfrac{nP + K_{pb}}{np + K_{pb}}\right]$ | $t_{1b21} = \bar{y}\left[\dfrac{nP - K_{pb}}{np - K_{pb}}\right]$ | 15.87 | 9.39 |
| $\beta_2(\phi)$ | P | $t_{1a22} = \bar{y}\left[\dfrac{\beta_2(\phi)P + P}{\beta_2(\phi)p + P}\right]$ | $t_{1b22} = \bar{y}\left[\dfrac{\beta_2(\phi)P - P}{\beta_2(\phi)p - P}\right]$ | 16.80 | 7.63 |
| NP | P | $t_{1a23} = \bar{y}\left[\dfrac{NP^2 + P}{NPp + P}\right]$ | $t_{1b23} = \bar{y}\left[\dfrac{NP^2 - P}{NPp - P}\right]$ | 15.93 | 9.26 |
| N | P | $t_{1a24} = \bar{y}\left[\dfrac{NP + P}{Np + P}\right]$ | $t_{1b24} = \bar{y}\left[\dfrac{NP - P}{Np - P}\right]$ | 13.23 | 11.32 |
| n | P | $t_{1a25} = \bar{y}\left[\dfrac{nP + P}{np + P}\right]$ | $t_{1b25} = \bar{y}\left[\dfrac{nP - P}{np - P}\right]$ | 14.51 | 10.28 |

**Appendix B.**

| | | | | | |
|---|---|---|---|---|---|
| **Some members (product-type) of the class** $t_1$ <br><br> **When** $w_0 = 0,\ w_1 = 1,\ w_2 = 0$ <br><br> $\alpha = -1,$ | | | | | |
| $K_1$ | $K_3$ | **Estimators** ($K_2 = 1$) | **Estimators** ($K_2 = -1$) | **PRE'S** <br><br> $K_2 = 1$ | **PRE'S** <br><br> $K_2 = -1$ |
| 1 | $C_p$ | $t_{1c1} = \overline{y}\left[\dfrac{p + C_p}{P + C_p}\right]$ | $t_{1d1} = \overline{y}\left[\dfrac{p - C_p}{P - C_p}\right]$ | 35.54 | 9.93 |
| 1 | $\beta_2(\phi)$ | $t_{1c2} = \overline{y}\left[\dfrac{p + \beta_2(\phi)}{P + \beta_2(\phi)}\right]$ | $t_{1d2} = \overline{y}\left[\dfrac{p - \beta_2(\phi)}{P - \beta_2(\phi)}\right]$ | 110.12 | 101.54 |
| $\beta_2(\phi)$ | $C_p$ | $t_{1c3} = \overline{y}\left[\dfrac{p\beta_2(\phi) + C_p}{P\beta_2(\phi) + C_p}\right]$ | $t_{1d3} = \overline{y}\left[\dfrac{p\beta_2(\phi) - C_p}{P\beta_2(\phi) - C_p}\right]$ | 6.09 | 0.127 |
| $C_p$ | $\beta_2(\phi)$ | $t_{1c4} = \overline{y}\left[\dfrac{pC_p + \beta_2(\phi)}{PC_p + \beta_2(\phi)}\right]$ | $t_{1d4} = \overline{y}\left[\dfrac{pC_p - \beta_2(\phi)}{PC_p - \beta_2(\phi)}\right]$ | 99.38 | 82.52 |
| 1 | $\rho_{pb}$ | $t_{1c5} = \overline{y}\left[\dfrac{p + \rho_{pb}}{P + \rho_{pb}}\right]$ | $t_{1d5} = \overline{y}\left[\dfrac{p - \rho_{pb}}{P - \rho_{pb}}\right]$ | 0.00135 | 5.42 |
| NP | $S_\phi$ | $t_{1c6} = \overline{y}\left[\dfrac{NPp + S_\phi}{NP^2 + S_\phi}\right]$ | $t_{1d6} = \overline{y}\left[\dfrac{NPp - S_\phi}{NP^2 - S_\phi}\right]$ | 2.53 | 1.23 |
| NP | f | $t_{1c7} = \overline{y}\left[\dfrac{NPp + f}{NP^2 + f}\right]$ | $t_{1d7} = \overline{y}\left[\dfrac{NPp - f}{NP^2 - f}\right]$ | 2.03 | 1.52 |
| $\beta_2(\phi)$ | $K_{pb}$ | $t_{1c8} = \overline{y}\left[\dfrac{\beta_2(\phi)p + K_{pb}}{\beta_2(\phi)P + K_{pb}}\right]$ | $t_{1d8} = \overline{y}\left[\dfrac{\beta_2(\phi)P - K_{pb}}{\beta_2(\phi)p - K_{pb}}\right]$ | 1.83 | 1.68 |
| NP | $K_{pb}$ | $t_{1c9} = \overline{y}\left[\dfrac{NPp + K_{pb}}{NP^2 + K_{pb}}\right]$ | $t_{1d9} = \overline{y}\left[\dfrac{NPp - K_{pb}}{NP^2 - K_{pb}}\right]$ | 1.89 | 1.63 |
| N | 1 | $t_{1c10} = \overline{y}\left[\dfrac{Np + 1}{NP + 1}\right]$ | $t_{1d10} = \overline{y}\left[\dfrac{Np - 1}{NP - 1}\right]$ | 2.37 | 1.30 |
| N | $C_p$ | $t_{1c11} = \overline{y}\left[\dfrac{Np + C_p}{NP + C_p}\right]$ | $t_{1d11} = \overline{y}\left[\dfrac{Np - C_p}{NP - C_p}\right]$ | 2.50 | 1.23 |
| N | $\rho_{pb}$ | $t_{1c12} = \overline{y}\left[\dfrac{Np + \rho_{pb}}{NP + \rho_{pb}}\right]$ | $t_{1d12} = \overline{y}\left[\dfrac{Np - \rho_{pb}}{NP - \rho_{pb}}\right]$ | 1.79 | 1.70 |
| N | $S_\phi$ | $t_{1c13} = \overline{y}\left[\dfrac{Np + S_\phi}{NP + S_\phi}\right]$ | $t_{1d13} = \overline{y}\left[\dfrac{Np - S_\phi}{NP - S_\phi}\right]$ | 2.16 | 1.44 |
| N | f | $t_{1c14} = \overline{y}\left[\dfrac{Np + f}{NP + f}\right]$ | $t_{1d14} = \overline{y}\left[\dfrac{Np - f}{NP - f}\right]$ | 1.98 | 1.56 |
| N | g=1-f | $t_{1c15} = \overline{y}\left[\dfrac{Np + g}{NP + g}\right]$ | $t_{1d15} = \overline{y}\left[\dfrac{Np - g}{NP - g}\right]$ | 2.34 | 1.32 |
| N | $K_{pb}$ | $t_{1c16} = \overline{y}\left[\dfrac{Np + K_{pb}}{NP + K_{pb}}\right]$ | $t_{1d16} = \overline{y}\left[\dfrac{Np - K_{pb}}{NP - K_{pb}}\right]$ | 1.93 | 1.60 |
| n | $\rho_{pb}$ | $t_{1c17} = \overline{y}\left[\dfrac{np + \rho_{pb}}{nP + \rho_{pb}}\right]$ | $t_{1d17} = \overline{y}\left[\dfrac{np - \rho_{pb}}{nP - \rho_{pb}}\right]$ | 1.49 | 1.96 |
| n | $S_\phi$ | $t_{1c18} = \overline{y}\left[\dfrac{np + \rho_{pb}}{nP + \rho_{pb}}\right]$ | $t_{1d18} = \overline{y}\left[\dfrac{np - \rho_{pb}}{nP - \rho_{pb}}\right]$ | 2.65 | 1.14 |
| n | f | $t_{1c19} = \bar{y}\left[\dfrac{np+f}{nP+f}\right]$ | $t_{1d19} = \bar{y}\left[\dfrac{np-f}{nP-f}\right]$ | 2.06 | 1.51 |
| n | g=1-f | $t_{1c20} = \bar{y}\left[\dfrac{np+g}{nP+g}\right]$ | $t_{1d20} = \bar{y}\left[\dfrac{np-g}{nP-g}\right]$ | 3.29 | 0.84 |
| n | $K_{pb}$ | $t_{1c21} = \bar{y}\left[\dfrac{np+K_{pb}}{nP+K_{pb}}\right]$ | $t_{1d21} = \bar{y}\left[\dfrac{np-K_{pb}}{nP-K_{pb}}\right]$ | 1.88 | 1.63 |
| $\beta_2(\phi)$ | P | $t_{1c22} = \bar{y}\left[\dfrac{\beta_2(\phi)p+P}{\beta_2(\phi)P+P}\right]$ | $t_{1d22} = \bar{y}\left[\dfrac{\beta_2(\phi)p-P}{\beta_2(\phi)P-P}\right]$ | 2.99 | 0.97 |
| NP | P | $t_{1c23} = \bar{y}\left[\dfrac{NPp+P}{NP^2+P}\right]$ | $t_{1d23} = \bar{y}\left[\dfrac{NPp-P}{NP^2-P}\right]$ | 2.37 | 1.30 |
| N | P | $t_{1c24} = \bar{y}\left[\dfrac{Np+P}{NP+P}\right]$ | $t_{1d24} = \bar{y}\left[\dfrac{Np-P}{NP-P}\right]$ | 2.11 | 1.47 |
| n | P | $t_{1c25} = \bar{y}\left[\dfrac{np+P}{nP+P}\right]$ | $t_{1d25} = \bar{y}\left[\dfrac{np-P}{nP-P}\right]$ | 2.49 | 1.23 |

**Appendix C.**

| **Some members (product-type) of the class** $t_2$ | | | |
|---|---|---|---|
| **When** $w_0 = 0,\ w_1 = 0,\ w_2 = 1$ | | | |

| $K_4$ | $K_5$ | **Estimators** $(\beta = 1, \lambda = -1)$ | PRE'S |
|---|---|---|---|
| 1 | $C_p$ | $t_{21} = \bar{y}\left(2 - \dfrac{p}{P}\exp\left[\dfrac{p-P}{p+P+2C_p}\right]\right)$ | 12.42 |
| 1 | $\beta_2(\phi)$ | $t_{22} = \bar{y}\left(2 - \dfrac{p}{P}\exp\left[\dfrac{p-P}{p+P+2\beta_2(\phi)}\right]\right)$ | 11.92 |
| $\beta_2(\phi)$ | $C_p$ | $t_{23} = \bar{y}\left(2 - \dfrac{p}{P}\exp\left[\dfrac{\beta_2(\phi)(p-P)}{\beta_2(\phi)(p+P)+2C_p}\right]\right)$ | 16.29 |
| $C_p$ | $\beta_2(\phi)$ | $t_{24} = \bar{y}\left(2 - \dfrac{p}{P}\exp\left[\dfrac{C_p(p-P)}{C_p(p+P)+2\beta_2(\phi)}\right]\right)$ | 12.53 |
| 1 | $\rho_{pb}$ | $t_{25} = \bar{y}\left(2 - \dfrac{p}{P}\exp\left[\dfrac{(p-P)}{(p+P)+2\rho_{pb}}\right]\right)$ | 13.86 |
| NP | $S_\phi$ | $t_{26} = \bar{y}\left(2 - \dfrac{p}{P}\exp\left[\dfrac{NP(p-P)}{NP(p+P)+2S_\phi}\right]\right)$ | 44.46 |
| NP | f | $t_{27} = \bar{y}\left(2 - \dfrac{p}{P}\exp\left[\dfrac{NP(p-P)}{NP(p+P)+2f}\right]\right)$ | 61.84 |

| $\beta_2(\phi)$ | $K_{pb}$ | $t_{28} = \overline{y}\left(2 - \dfrac{p}{P}\exp\left[\dfrac{\beta_2(\phi)(p-P)}{\beta_2(\phi)(p+P)+2K_{pb}}\right]\right)$ | 40.17 |
|---|---|---|---|
| $NP$ | $K_{pb}$ | $t_{29} = \overline{y}\left(2 - \dfrac{p}{P}\exp\left[\dfrac{NP(p-P)}{NP(p+P)+2K_{pb}}\right]\right)$ | 48.09 |
| $N$ | $1$ | $t_{210} = \overline{y}\left(2 - \dfrac{p}{P}\exp\left[\dfrac{N(p-P)}{N(p+P)+2}\right]\right)$ | 54.10 |
| $N$ | $C_p$ | $t_{211} = \overline{y}\left(2 - \dfrac{p}{P}\exp\left[\dfrac{N(p-P)}{N(p+P)+2C_p}\right]\right)$ | 44.84 |
| $N$ | $\rho_{pb}$ | $t_{212} = \overline{y}\left(2 - \dfrac{p}{P}\exp\left[\dfrac{N(p-P)}{N(p+P)+2\rho_{pb}}\right]\right)$ | 56.48 |
| $N$ | $S_\phi$ | $t_{213} = \overline{y}\left(2 - \dfrac{p}{P}\exp\left[\dfrac{N(p-P)}{N(p+P)+2S_\phi}\right]\right)$ | 62.40 |
| $N$ | $f$ | $t_{214} = \overline{y}\left(2 - \dfrac{p}{P}\exp\left[\dfrac{N(p-P)}{N(p+P)+2f}\right]\right)$ | 65.67 |
| $N$ | $g=1-f$ | $t_{215} = \overline{y}\left(2 - \dfrac{p}{P}\exp\left[\dfrac{N(p-P)}{N(p+P)+2g}\right]\right)$ | 54.47 |
| $N$ | $K_{pb}$ | $t_{216} = \overline{y}\left(2 - \dfrac{p}{P}\exp\left[\dfrac{N(p-P)}{N(p+P)+2K_{pb}}\right]\right)$ | 63.21 |
| $N$ | $\rho_{pb}$ | $t_{217} = \overline{y}\left(2 - \dfrac{p}{P}\exp\left[\dfrac{n(p-P)}{n(p+P)+2\rho_{pb}}\right]\right)$ | 38.59 |

| N | $S_\phi$ | $t_{218} = \overline{y}\left(2 - \dfrac{p}{P}\exp\left[\dfrac{n(p-P)}{n(p+P)+2S_\phi}\right]\right)$ | 52.19 |
|---|---|---|---|
| N | f | $t_{219} = \overline{y}\left(2 - \dfrac{p}{P}\exp\left[\dfrac{n(p-P)}{n(p+P)+2f}\right]\right)$ | 63.74 |
| N | g=1-f | $t_{220} = \overline{y}\left(2 - \dfrac{p}{P}\exp\left[\dfrac{n(p-P)}{n(p+P)+2g}\right]\right)$ | 35.33 |
| N | $K_{pb}$ | $t_{221} = \overline{y}\left(2 - \dfrac{p}{P}\exp\left[\dfrac{n(p-P)}{n(p+P)+2K_{pb}}\right]\right)$ | 54.68 |
| $\beta_2(\phi)$ | P | $t_{222} = \overline{y}\left(2 - \dfrac{p}{P}\exp\left[\dfrac{\beta_2(\phi)(p-P)}{\beta_2(\phi)(p+P)+2P}\right]\right)$ | 47.53 |
| NP | P | $t_{223} = \overline{y}\left(2 - \dfrac{p}{P}\exp\left[\dfrac{NP(p-P)}{NP(p+P)+2P}\right]\right)$ | 54.10 |
| N | P | $t_{224} = \overline{y}\left(2 - \dfrac{p}{P}\exp\left[\dfrac{N(p-P)}{N(p+P)+2P}\right]\right)$ | 64.43 |
| N | P | $t_{225} = \overline{y}\left(2 - \dfrac{p}{P}\exp\left[\dfrac{n(p-P)}{n(p+P)+2P}\right]\right)$ | 58.91 |

In addition to above estimators a large number of estimators can also be generated from the proposed estimators just by putting different values of constants $w_i$'s, $K_1$, $K_2$, $K_3$, $K_4$, $K_5$, $\alpha$, $\beta$ and $\lambda$.

# Two-Warehouse Fuzzy Inventory Model
# with K-Release Rule

**Neeraj Kumar\*, Sanjey Kumar\* and Florentin Smarandache\*\***

*Department of Mathematics, SRM University

Delhi NCR, Sonepat, Haryana

**Department of Mathematics, University of New Mexico

Gallup, NM 87301, USA

*Corresponding Author, sanjeysrm1984@gmail.com

**ABSTRACT:** Fuzzy set theory is primarily concerned with how to quantitatively deal with imprecision and uncertainty, and offers the decision maker another tool in addition to the classical deterministic and probabilistic mathematical tools that are used in modeling real-world problems. The present study investigates a fuzzy economic order quantity model for two storage facility. The demand, holding cost, ordering cost, storage capacity of the own - warehouse are taken as a trapezoidal fuzzy numbers. Graded Mean Representation is used to defuzzify the total cost function and the results obtained by this method are compared with the help of a numerical example. Sensitivity analysis is also carried out to explore the effect of changes in the values of some of the system parameters. The proposed methodology is applicable to other inventory models under uncertainty.

**Keywords:** Inventory, Two – warehouse system, Fuzzy Variable, Trapezoidal Fuzzy Number, Graded mean representation method and K – release rule.

## 1. INTRODUCTION

In most of the inventory models that had been proposed in the early literature, the associated costs are assumed to be precise, although the real-world inventory costs usually exist with imprecise components. In this case, customer demand as one of the key parameters and source of uncertainty have been most often treated by a probability distribution. However, the probability-based approaches may not be sufficient enough to reflect all uncertainties that may arise in a real-world inventory system. Modelers may face some difficulties while trying to build a valid model of an inventory system, in which the related costs cannot be determined precisely. For example, costs may be dependent on some foreign monetary unit. In such a case, due to a change in the exchange rates, the costs are often not known precisely.

Fuzzy set theory, originally introduced by Zadeh [1], provides a framework for considering parameters that are vaguely or unclearly defined or whose values are imprecise or determined based on subjective beliefs of individuals. Petrovic et al. [2] presented newsboy problem assuming that demand and backorder cost are fuzzy numbers. Kaufamann and Gupta [3] introduced to fuzzy arithmetic: theory and application. The application of fuzzy theory to inventory problem has been proposed by Kacprzyk and Staniewski [4]. Roy and Maiti [5] presented a fuzzy inventory model with constraint. Roy and Maiti [6] developed a fuzzy EOQ model with demand-dependent unit cost under limited storage capacity. Ishii and Konno [7] introduced fuzziness of shortage cost explicitly into classical inventory problem. Chen and Hsieh [8] established a fuzzy economic production model to treat the inventory problem with all the parameters and variables, which are fuzzy numbers. Hsieh [9] presented a fuzzy production inventory model. Yao and Chiang [10] presented an inventory model without backorder with fuzzy total cost and fuzzy storing cost defuzzified by centroid and signed distance. Dutta et al. [11] developed a single-period inventory model with fuzzy random variable demand. In that study, they have applied graded mean integration representation method to find the optimum order quantity. Chen and Chang [12] presented an optimization of fuzzy production inventory model. In this study, they have used 'Function Principle' as arithmetical operations of fuzzy total production inventory cost and also used the 'Graded Mean Integration Representation method' to defuzzify the fuzzy total production and inventory cost. Mahata and Goswami [13] presented a fuzzy inventory model for deteriorating items with the help of fuzzy numbers and so on.

Most of the classical inventory models discussed in the literature deals with the situation of a single warehouse. Because of capacity limitation a single warehouse would not be always sufficient. Additional warehouse are necessary to store excess items. Therefore due to the limited capacity of the existing warehouse (Rented warehouse, RW) is acquired to keep excess items. In practice, large stock attracts the management due to either an attractive price discount for bulk purchase or the acquisition cost being higher than the holding cost in RW. The actual service to the customer is done at OW only. Usually the holding cost is greater in RW than in OW. So in order to reduce the holding cost. The stock of rented warehouse is transferred to the own warehouse. Hartley [14] was discussed a model under the assumption that the cost of transporting a unit from RW to OW is not significantly high. It was as the case with two levels of storage. Sarma [15] extended the model with two levels of storage given by Hartley, by considering the transportation cost of a unit from rented warehouse to own warehouse. Maurdeswar and Sathe [16] discussed this model by relaxing the condition on production rate (finite production rate). Dave [17] considered it for finite and infinite replenishment, assuming the cost of transportation depending on the quantity to be transported. Pakkala and Achary [18] developed a model for deteriorating items with two warehouses. They extended it with bulk release rule, after words, Gowsami and Chaudhari [19] formulated models for time dependent demand. Kar et al. [20] suggested a two level inventory model for linear trend in demand. Yang [21] considered a two-warehouse inventory models for deteriorating items with shortages under inflation. Singh et al. [22] presented two-warehouse inventory model without shortage for exponential demand rate and an optimum release rule. Jaggi and Verma [23] developed a deterministic order level inventory model with two storage facilities. It has been observed in supermarkets that the demand rate is usually influenced by the amount of stock level, that is, the demand rate may go up or down with the on-hand stock level. Singh et al. [24] developed a deterministic two-warehouse inventory model for deteriorating items with stock-dependent demand and shortages. Neeraj et al. [25] developed three echelon supply chain inventory model with two storage facility. Neeraj et al. [26] presented a two-warehouse inventory model with K-release rule and learning effect. Neeraj et al. [27] considered effect of salvage value on a two-warehouse inventory model. Recently, Kumar and Kumar [28] developed an inventory model with stock dependent demand rate for deterioration items.

Here, in this paper the cost of transporting a unit is considered to be significant and the effect of releasing the stocks of RW in n shipments with a bulk size of K units per shipment, instead of withdrawing an arbitrary quantity, is assumed. Here, K is to be decided optimally and is call this as K-release rule. This problem is to decide the optimal values of Q and C, which minimize the sum of ordering, holding and transportation costs of the system. Here, we assumed that the storage capacity of the own – warehouse, the holding cost in both warehouses and ordering cost is fuzzy in nature. The associated total cost minimization is illustrated by numerical example and sensitivity analysis is carried out by using *MATHEMATICA–5.2* for the feasibility and applicability of our model.

## 2. ASSUMPTIONS AND NOTATIONS:

The following assumptions are used to analyze this inventory model:

1.  D is the constant demand rate.
2.  W is the storage capacity of the OW.

3.  A is the fixed set – up cost per order.

4.  C(Q) is the cost function.

5.  Q is the highest inventory level.

6.  H is the holding cost in OW.

7.  F is the holding cost in RW.

8.  $\tilde{D}$ is the fuzzy demand rate.

9.  $\tilde{A}$ is the fuzzy set – up cost per order.

10. $\tilde{H}$ is the fuzzy holding cost in OW.

11. $\tilde{F}$ is the fuzzy holding cost in RW.

12. $\tilde{C}(Q)$ is the fuzzy cost function.

13. $\tilde{W}$ is the fuzzy storage capacity of the OW.

14. The holding cost per unit in OW is higher than in RW.

15. The storage capacity of OW as W and that of RW is unlimited.

16. The transportation cost of K units from RW to OW is $C_t$ at a time, which is constant over time.

17. The items of RW are transferred to OW in 'n' shipments of which K (K ≤ W) units are transported in each shipment.

18. Replenishment rate is infinite.

19. Lead-time is zero.

20. Consumption takes place only in OW.

## 3. FUZZY SETS, MEMBERSHIP FUNCTION, DEFUZZIFYING APPROACH AND ARITHMETICAL OPERATIONS

### 3.1. Fuzzy Sets

A fuzzy set is a class of objects with a continuum of grades of membership. Such a set is characterized by a membership (characteristic) function which assigns to each object a grade of membership ranging between zero and one. Let X={x} denote a space of objects. Then a fuzzy set $A$ in $X$ is a set of ordered pairs:

$$A = \left\{ x, \mu_A(x) \right\}, \; x \in X$$

Where, $\mu_A(x)$ is termed " the grade of the membership of $x$ in $A$ ". For simplicity, $\mu_A(x)$ is a number in the interval [0, 1], with the grades of unity and zero respectively, full membership and non-membership in the fuzzy set. An object (point) P contained in a set (class) Q is an element of $Q$ $(P \subset Q)$.

### 3.2. Membership Function

Membership Function

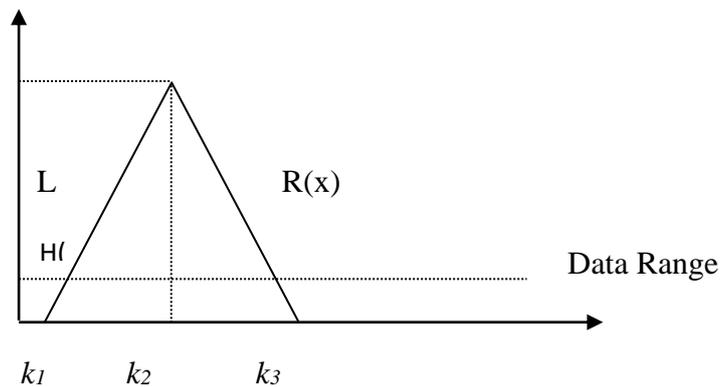

**Fig. 1 Membership function for triangle number**

At the outset it would be prudent introduce the concept of membership function. There are different shapes of membership function in the inventory control such as the triangle and trapezoid. The shapes of the triangle membership function and the trapezoid membership function are shown in Fig. 1 and 2.

Membership Function

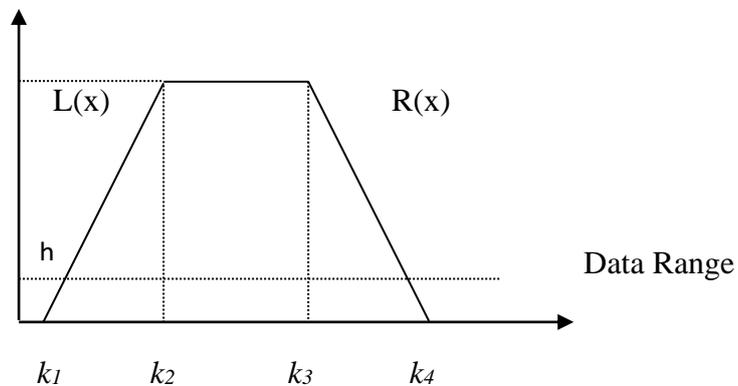

**Fig. 2 Membership function for trapezoid number**

Ã is assumed as a fuzzy number. If Ã is a triangle number, Ã can be represented as Ã = $[k_1, k_2, k_3]$ subject to the constraint $0 < k_1 \leq k_2 \leq k_3$. While Ã is a trapezoid fuzzy number, Ã = $[k_1, k_2, k_3, k_4]$ subject to the constraint that $0 < k_1 \leq k_2 \leq k_3 \leq k_4$. Membership function of the triangle and trapezoid fuzzy numbers can be defined as follows:

$$\mu_{\tilde{A}}(x) = \begin{cases} 0 & x < k_1, \ x > k_3 \\ L(x) = \dfrac{x - k_1}{k_2 - k_1} & k_1 \leq x < k_2 \\ R(x) = \dfrac{k_3 - x}{k_3 - k_2} & k_2 \leq x < k_3 \end{cases}$$

$$\mu_{\tilde{A}}(x) = \begin{cases} 0 & x < k_1, x > k_4 \\ L(x) = \dfrac{x - k_1}{k_2 - k_1} & k_1 \leq x < k_2 \\ 1 & k_2 \leq x < k_3 \\ R(x) = \dfrac{k_4 - x}{k_4 - k_3} & k_3 \leq x \leq k_4 \end{cases}$$

where $\mu_{\tilde{A}}(x)$ is a membership function.

### 3.3. Graded Mean Integration Representation Method

In this study, generalized fuzzy number Ã was denoted in Fig. 6.1 as Ã = $\left( c, a, b, d, \omega_A \right)_{LR}$. When $\omega_A = 1$, we simplify the notation as $\tilde{A} = \left( c, a, b, d \right)_{LR}$. Chen and Hsieh (1999) introduced the graded mean integration representation method of generalized fuzzy number based on the integral value of graded mean $h$–level of generalized fuzzy number. Its meaning is as follows:

Let $L^{-1}$ and $R^{-1}$ are inverse function of $L$ and $R$ respectively, then the graded mean $h$–level value of generalized fuzzy number $\tilde{A} = \left( c, a, b, d, W_A \right)_{LR}$ is $h \left( L^{-1}(h) + R^{-1}(h) \right) / 2$ as Fig. 3.

Then the graded mean integration representation of Ã is

$$P(\tilde{A}) = \int_o^{W_A} \frac{h \left( L^{-1}(A) + R^{-1}(h) \right)}{2} dh \left/ \int_o^{W_A} h \, dh \right. ,$$

where $0 < h \leq W_A$ and $0 < W_A \leq 1$.

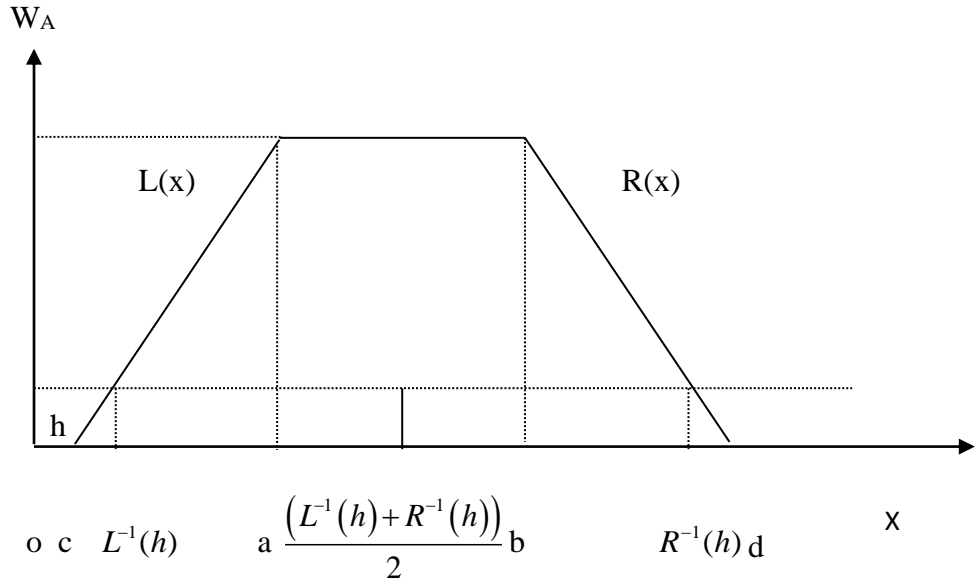

**Fig. 3 The graded mean h-level of generalized fuzzy number A = (c, a, b, d, W$_A$)$_{LR}$**

In the present, the generalized trapezoidal fuzzy number has been used as the type of all fuzzy parameters in our proposed inventory models. The very popular generalized trapezoidal fuzzy number $B$ is a special case of generalized fuzzy number and can be denoted as $\tilde{B} = \left(c, a, b, d; W_B\right)$ its' corresponding graded mean integration representation is

$$P(\tilde{B}) = \int_0^{W_B} \frac{h(c+d+\left(a-c-d+b\right)h/W_B)dh}{2} \Bigg/ \int_0^{W_B} hdh = \frac{c+2a+2b+d}{6}$$

where $a, b, c, d$ are any real numbers.

### 3.4. Properties of Second Function Principle

Chen (1985) proposed second function principal to be as the fuzzy arithmetical operations between generalized trapezoidal fuzzy numbers. Because it does not change the type of membership function of generalized fuzzy number after arithmetical operations. It reduces the trouble and tediousness of operations. Furthermore, Chen already proved the properties of fuzzy arithmetical operations under second function principle. Here some properties of the fuzzy arithmetical operations have been described as follows:

Suppose $\tilde{A}_1 = (c_1, a_1, b_1, d_1)$ and $\tilde{A}_2 = (c_2, a_2, b_2, d_2)$ are two generalized trapezoidal fuzzy numbers. Then

1. The addition of $\tilde{A}_1$ and $\tilde{A}_2$ is $\tilde{A}_1 \oplus \tilde{A}_2 = (c_1 + c_2, a_1 + a_2, b_1 + b_2, d_1 + d_2)$

2. The multiplication of $\tilde{A}_1$ and $\tilde{A}_2$ is $\tilde{A}_1 \otimes \tilde{A}_2 = (c_1 c_2, a_1 a_2, b_1 b_2, d_1 d_2)$

3. $-\tilde{A}_2 = (-d_2, -b_2, -a_2, -c_2)$ Then the subtraction of $\tilde{A}_1$ and $\tilde{A}_2$ is $\tilde{A}_1 \ominus \tilde{A}_2 = (c_1 - d_2, a_1 - b_2, b_1 - a_2, d_1 - c_2)$

4. $1/\tilde{A}_2 = A^{-1}{}_2 = \left( \dfrac{1}{d_2}, \dfrac{1}{b_2}, \dfrac{1}{a_2}, \dfrac{1}{c_2} \right)$ where $c_2, a_2, b_2$ and $d_2$ are all positive real numbers. If $c_1, a_1, b_1, d_1, c_2, a_2, b_2$ and $d_2$ are all non zero positive real numbers, then the division of $\tilde{A}_1$ and $\tilde{A}_2$ is $\tilde{A}_1 \oslash \tilde{A}_2 = \left( \dfrac{c_1}{d_2}, \dfrac{a_1}{b_2}, \dfrac{b_1}{a_2}, \dfrac{d_1}{c_2} \right)$.

## 4. MODEL DEVELOPMENT

Initially the company ordered Q units of the item, out of which W units is kept in OW and Z units are kept in RW, where Z = (Q - W). Initially, demand is satisfied using the stocks of OW until the stock level drops to (W-K) units. At this stage, K units from RW are transported to OW to meet further demand and this process is repeated 'n' times until the stocks of RW are exhausted. The remaining (W-K) units in OW are used again at this stage. The inventory situation in RW and OW are shown in the figure 1.

The inventory units in RW can be seen to be equal.

$$A_t = t_{ik} \left[ Z + (Z-K) + (Z-2K) + \ldots\ldots + (Z-(n-1)K) \right] = t_{ik} \frac{Z(n+1)}{2} \tag{4.1}$$

Where $t_{ik}$ = K/D, the time taken for the consumption of K units, since Z = (Q - W) and the holding cost in RW is F(i), we have-

$$FA_t = Ft_{ik} \frac{Z(n+1)}{2} = \frac{FK}{D}(Q-W)\frac{(n+1)}{2} = \frac{FK(n+1)(Q-W)}{2D} \tag{4.2}$$

The cost of transporting the units from RW to OW in 'n' shipments is given by

$$nC_t = (Z/K)C_t \qquad (4.3)$$

Since $n = Z/K$

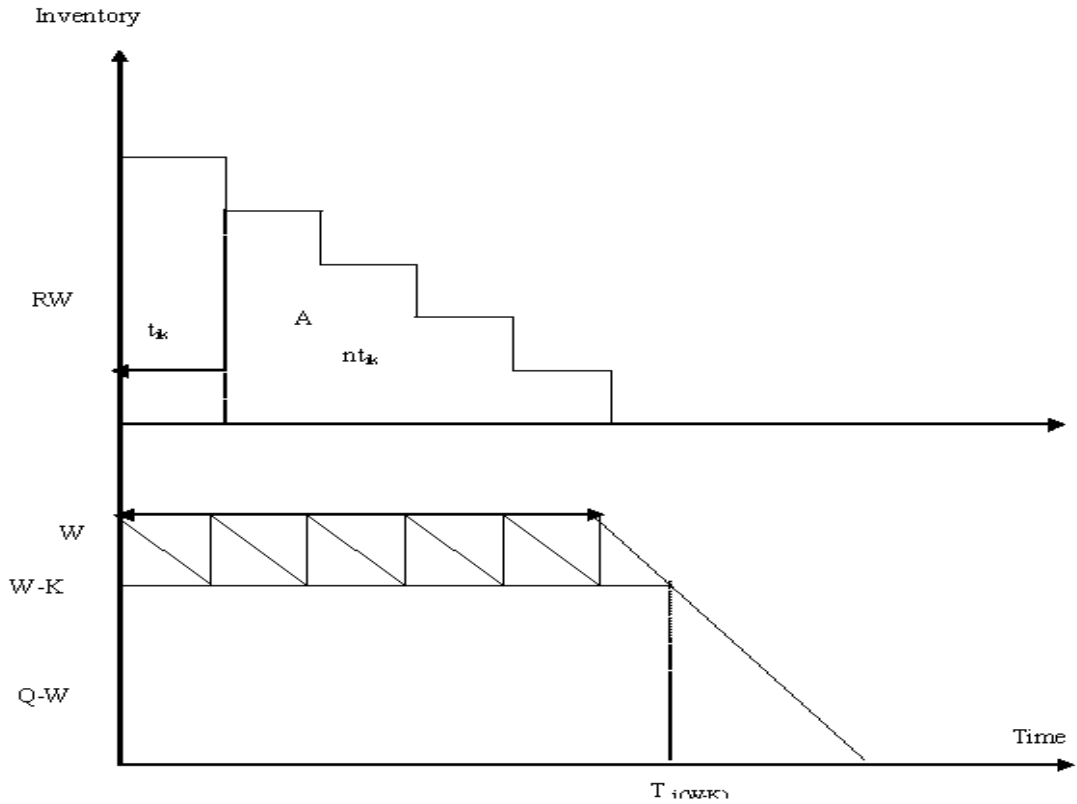

**Fig. 1: Graphical representation of two-storage inventory model**

When K units are drawn from RW in each shipment, more are carried in OW for a period of $t_k$ and hence account for a holding cost of KH (i) $t_{ki}$ / 2. Since there are 'n' such shipments and taking into consideration, the initial K units of OW, the holding cost for these items is

$$(n+1)HKt_{ik}/2 = (n+1)HK^2/2D \qquad (4.4)$$

A quantity of (W- K) units is kept unused in OW for a period of $t_{i(W-K)} = (n+1)t_{ik}$ and an average inventory during usage in OW is (W - K)/2 units for a period $(t - t_{i(W-K)})$. Hence the inventory holding cost in OW for these items is

$$H[K(W-K)(n+1)/D + (W - K)^2/2D]. \qquad (4.5)$$

The fixed ordering cost per order is A. Then the total inventory cost for the system using (4.2) to (4.5) becomes

$$C = A + (n+1)\frac{FK(Q-W)}{2D} + \frac{K^2 H}{2D} + \frac{HK(W-K)}{2D} + nC_t + \frac{(W-K)^2 H}{2D} \qquad (4.6)$$

The average inventory cost

C(Q, K) = C / t

But we have t = Q / D, Z = Q – W and n = Z / K = Q – W/K

Total average cost becomes

$$C(Q,K) = \frac{AD}{Q} + \frac{FQ}{2} - W(F-H) + \frac{K}{2}(F-H) - \frac{KW}{2Q}(F-H) + C_t\frac{(Q-W)D}{QK} + \frac{W^2}{2Q}(F-H) \qquad (4.7)$$

**Fuzzy Model:** Due to uncertainly in the environment it is not easy to define all the parameters precisely, accordingly we assume some of these parameters namely $\tilde{D}$, $\tilde{F}$, $\tilde{H}$, $\tilde{A}$ and $\tilde{W}$ may change within some limits. Let $\tilde{D} = [d_1, d_2, d_3, d_4]$, $\tilde{F} = [f_1, f_2, f_3, f_4]$, $\tilde{H} = [h_1, h_2, h_3, h_4]$, $\tilde{A} = [a_1, a_2, a_3, a_4]$, $\tilde{W} = [w_1, w_2, w_3, w_4]$ are as trapezoidal fuzzy numbers. In this case, the total fuzzy cost per unit time is given by

$$\tilde{C}(Q,K) = \left( \left( \tilde{A} \otimes \tilde{D} \right) \oslash Q \right) \oplus \left( \left( \tilde{F} \otimes Q \right) \oslash 2 \right) \ominus \left( \tilde{W} \otimes \left( \tilde{F} \ominus \tilde{H} \right) \right) \oplus \left( \left( K \otimes \left( \tilde{F} \ominus \tilde{H} \right) \right) \oslash 2 \right) \ominus \left( \left( K \otimes \left( \tilde{W} \otimes \left( \tilde{F} \ominus \tilde{H} \right) \right) \right) \oslash 2Q \right)$$
$$\oplus \left( \left( \left( C_t \otimes \left( Q \ominus \tilde{W} \right) \right) \otimes \tilde{D} \right) \oslash QK \right) \oplus \left( \left( \tilde{W} \otimes \left( \tilde{W} \otimes \left( \tilde{F} \ominus \tilde{H} \right) \right) \right) \oslash 2Q \right) \qquad (4.8)$$

By second function principal, one has

$$\tilde{C}(Q,K) = \left( \frac{a_1 d_1}{Q} + \frac{f_1 Q}{2} - (f_4 - h_1)w_4 + \frac{K(f_1 - h_4)}{2} - \frac{K(f_4 - h_1)w_4}{2Q} + \frac{C_t(Q - w_4)d_1}{QK} + \frac{(f_1 - h_4)w_1^2}{2Q}, \right.$$

$$\frac{a_2 d_2}{Q} + \frac{f_2 Q}{2} - (f_3 - h_2)w_3 + \frac{K(f_2 - h_3)}{2} - \frac{K(f_3 - h_2)w_3}{2Q} + \frac{C_t(Q - w_3)d_2}{QK} + \frac{(f_2 - h_3)w_1^2}{2Q},$$

$$\frac{a_3 d_3}{Q} + \frac{f_3 Q}{2} - (f_2 - h_3)w_2 + \frac{K(f_3 - h_2)}{2} - \frac{K(f_2 - h_3)w_2}{2Q} + \frac{C_t(Q - w_2)d_3}{QK} + \frac{(f_3 - h_2)w_1^2}{2Q},$$

$$\left. \frac{a_4 d_4}{Q} + \frac{f_4 Q}{2} - (f_1 - h_4)w_1 + \frac{K(f_4 - h_1)}{2} - \frac{K(f_1 - h_4)w_1}{2Q} + \frac{C_t(Q - w_1)d_4}{QK} + \frac{(f_4 - h_1)w_1^2}{2Q} \right)$$

Now we defuzzify the total cost per unit time, using graded mean integration representation method, the result is

$$P\left(\tilde{C}(Q,K)\right)=\frac{1}{6}\left[\left(\frac{a_1 d_1}{Q}+\frac{f_1 Q}{2}-(f_4-h_1)w_4+\frac{K(f_1-h_4)}{2}-\frac{K(f_4-h_1)w_4}{2Q}+\frac{C_t(Q-w_4)d_1}{QK}+\frac{(f_1-h_4)w_1^2}{2Q}\right.\right.$$

$$+2\left(\frac{a_2 d_2}{Q}+\frac{f_2 Q}{2}-(f_3-h_2)w_3+\frac{K(f_2-h_3)}{2}-\frac{K(f_3-h_2)w_3}{2Q}+\frac{C_t(Q-w_3)d_2}{QK}+\frac{(f_2-h_3)w_1^2}{2Q}\right)$$

$$+2\left(\frac{a_3 d_3}{Q}+\frac{f_3 Q}{2}-(f_2-h_3)w_2+\frac{K(f_3-h_2)}{2}-\frac{K(f_2-h_3)w_2}{2Q}+\frac{C_t(Q-w_2)d_3}{QK}+\frac{(f_3-h_2)w_1^2}{2Q}\right)$$

$$\left.\left.+\frac{a_4 d_4}{Q}+\frac{f_4 Q}{2}-(f_1-h_4)w_1+\frac{K(f_4-h_1)}{2}-\frac{K(f_1-h_1)w_1}{2Q}+\frac{C_t(Q-w_1)d_4}{QK}+\frac{(f_4-h_1)w_1^2}{2Q}\right)\right]$$

(4.9)

The optimal values of Q and K, which minimizes (4.8), are obtained by solving

$$\frac{\partial P\left(\tilde{C}(Q,K)\right)}{\partial Q}=0 \ and \ \frac{\partial P\left(\tilde{C}(Q,K)\right)}{\partial K}=0$$

(4.10)

we get

$$Q=\left[\frac{\begin{array}{l}2a_1 d_1+4a_2 d_2+4a_3 d_3+2a_4 d_4-K\left[(f_4-h_1)w_4+(f_3-h_2)w_3+(f_2-h_3)w_2+(f_1-h_4)w_1\right]\\-\frac{2C_t}{K}\left(w_4 d_1+w_3 d_2+w_2 d_3+w_1 d_4\right)+(f_4-h_1)w_4^2+(f_3-h_2)w_3^2+(f_2-h_3)w_2^2+(f_1-h_4)w_1^2\end{array}}{f_1+f_2+f_3+f_4}\right]^{1/2}$$

(4.11)

and

$$K=\sqrt{\frac{\frac{C_t}{Q}\left[(Q-w_4)+2(Q-w_3)+2(Q-w_2)+(Q-w_1)\right]}{\frac{(f_1-h_4)}{2}+(f_2-h_3)+(f_3-h_2)+\frac{(f_4-h_1)}{2}-\frac{(f_1-h_4)w_1}{2Q}-\frac{(f_2-h_3)w_2}{Q}-\frac{(f_3-h_2)w_3}{Q}-\frac{(f_4-h_1)w_4}{2Q}}}$$

(4.12)

## 5. COST-REDUCTION DUE TO K-RELEASE RULE

The unit cost of transportation with K-release rule is $C_t^{'} = C_t / K$ . Suppose the unit cost of transportation is $C_t^{*}$ without bulk transportation. The bulk transportation will be economical only if $C_t^{*} > C_t^{'}$ . Hence without K-release rule, the cost function becomes-

$$C(Q) = \frac{AD}{Q} + \frac{FQ}{2} + \frac{W^2(F-H)}{2Q} - W(F-H) + \frac{(Q-W)C_t^{*}D}{Q} \tag{5.1}$$

**Fuzzy Model:** Due to uncertainly in the environment it is not easy to define all the parameters precisely, accordingly we assume some of these parameters namely $\tilde{D}$, $\tilde{F}$, $\tilde{H}$, $\tilde{A}$ and $\tilde{W}$ may change within some limits.

Let $\tilde{D} = [d_1, d_2, d_3, d_4]$, $\tilde{F} = [f_1, f_2, f_3, f_4]$, $\tilde{H} = [h_1, h_2, h_3, h_4]$, $\tilde{A} = [a_1, a_2, a_3, a_4]$, $\tilde{W} = [w_1, w_2, w_3, w_4]$ are as trapezoidal fuzzy numbers. In this case, the total fuzzy cost per unit time is given by

$$\tilde{C}(Q) = \left( (\tilde{A} \otimes \tilde{D}) \oslash Q \right) \oplus \left( (\tilde{F} \otimes Q) \oslash 2 \right) \oplus \left( \left( \tilde{W} \otimes (\tilde{W} \otimes (\tilde{F} \ominus \tilde{H})) \right) \right) \oslash 2Q \right) \ominus \left( \tilde{W} \otimes (\tilde{F} \ominus \tilde{H}) \right)$$
$$\oplus \left( \left( \left( (Q \ominus \tilde{W}) \otimes C_t^{*} \right) \otimes \tilde{D} \right) \oslash Q \right) \tag{5.2}$$

By second function principal, one has

$$\tilde{C}(Q) = \left( \frac{a_1 d_1}{Q} + \frac{f_1 Q}{2} - (f_4 - h_1)w_4 + \frac{C_t^{*}(Q - w_4)d_1}{Q} + \frac{(f_1 - h_4)w_1^2}{2Q} , \right.$$

$$\frac{a_2 d_2}{Q} + \frac{f_2 Q}{2} - (f_3 - h_2)w_3 + \frac{C_t^{*}(Q - w_3)d_2}{Q} + \frac{(f_2 - h_3)w_1^2}{2Q} ,$$

$$\frac{a_3 d_3}{Q} + \frac{f_3 Q}{2} - (f_2 - h_3)w_2 + \frac{C_t^{*}(Q - w_2)d_3}{Q} + \frac{(f_3 - h_2)w_1^2}{2Q} ,$$

$$\left. \frac{a_4 d_4}{Q} + \frac{f_4 Q}{2} - (f_1 - h_4)w_1 + \frac{C_t^{*}(Q - w_1)d_4}{Q} + \frac{(f_4 - h_1)w_1^2}{2Q} \right)$$

Now we defuzzify the total cost per unit time, using graded mean integration representation method, the result is

$$P\left(\tilde{C}(Q)\right) = \frac{1}{6} \left\{ \begin{array}{l} \left( \dfrac{a_1 d_1}{Q} + \dfrac{f_1 Q}{2} - (f_4 - h_1) w_4 + \dfrac{C_t^*(Q - w_4) d_1}{Q} + \dfrac{(f_1 - h_4) w_1^2}{2Q} \right) \\[2mm] + 2 \left( \dfrac{a_2 d_2}{Q} + \dfrac{f_2 Q}{2} - (f_3 - h_2) w_3 + \dfrac{C_t^*(Q - w_3) d_2}{Q} + \dfrac{(f_2 - h_3) w_1^2}{2Q} \right) \\[2mm] + 2 \left( \dfrac{a_3 d_3}{Q} + \dfrac{f_3 Q}{2} - (f_2 - h_3) w_2 + \dfrac{C_t^*(Q - w_2) d_3}{Q} + \dfrac{(f_3 - h_2) w_1^2}{2Q} \right) \\[2mm] + \left( \dfrac{a_4 d_4}{Q} + \dfrac{f_4 Q}{2} - (f_1 - h_4) w_1 + \dfrac{C_t^*(Q - w_1) d_4}{Q} + \dfrac{(f_4 - h_1) w_1^2}{2Q} \right) \end{array} \right\} \tag{5.3}$$

The optimal value of Q, which minimizes (5.1), is obtained by $\dfrac{dP\left(\tilde{C}(Q)\right)}{dQ} = 0$

$$Q = \left[ \frac{-a_1 d_1 - 2a_2 d_2 - 2a_3 d_3 - a_4 d_4 + C_t \left( w_4 d_1 + 2w_3 d_2 + 2w_2 d_3 + w_1 d_4 \right)}{-\dfrac{(f_4 - h_1)}{2} w_4{}^2 - (f_3 - h_2) w_3{}^2 - (f_2 - h_3) w_2{}^2 - \dfrac{(f_1 - h_4)}{2} w_1{}^2} \right]^{1/2} \tag{5.4}$$

The proposed K-released rule will be economical if

$$\left[ C(Q) - C(Q, K) \right] > 0$$

From equation (4.7) and (4.8) we see that-

$$\left[ C(Q) - C(Q, K) \right] = \left( 1 - \frac{W}{Q} \right) \left[ D \left( C_t^* - C_t^{'} \right) - \frac{K}{2} (F - H) \right] \tag{5.5}$$

and hence the inequality

$$\left( 1 - \frac{W}{Q} \right) \left[ D \left( C_t^* - C_t^{'} \right) - \frac{K}{2} (F - H) \right] > 0$$

$$\Rightarrow \left( C_t^* - C_t^{'} \right) > \frac{K(F - H)}{2D} \tag{5.6}$$

must be satisfied.

Thus for a given situation, if the unit cost of transportation with bulk release rule satisfies the inequality (5.6), K-release rule must be economical.

## 6. NUMERICAL EXAMPLE

Consider an inventory system with following parametric values:

**Crisp Model:** demand rate D = 2000, $C_t$ = 0.5, F = 8.5, H = 7.5, W = 100, A = 150. With the help of the above values, we find the optimal values of ordering quantity and total cost with and without K- release which is given as:

With K – release rule: Q = 221.62 & C (Q, K) = 3456.46

And without K – release rule: Q = 216.68 & C (Q, K) = 3585.43

**Fuzzy Model:** $\tilde{D}$ = [1900, 2000, 2000, 1900], $\tilde{F}$ = [8.075, 8.5, 8.5, 8.075], $\tilde{H}$ = [7.125, 7.5, 7.5, 7.125], $\tilde{A}$ = [142.5, 150, 150, 142.5], $\tilde{W}$ = [95, 100, 100, 95]. The optimal values of ordering quantity and total cost with and without K- release which is given as:

With K – release rule: Q = 225.62 & C (Q, K) = 3458.46

And without K – release rule: Q = 210.68 & C (Q, K) = 3587.43

## 7. CONCLUSION

Two storage inventory models discussed in this paper and developed under the assumption that the distribution of the items to the customers takes place at OW only. Because of the distance factor, it is natural to consider the transportation cost associated with the transfer of items from RW to OW. Further, the concept of K-release rule is more pragmatic, as holding large inventory in RW is every expensive. With the help of numerical examples, it is clear that the effect of fuzzy cannot be ignored. We can earn more profit by consider the effect of fuzzy on ordering and holding cost in each lot. This model gives the direction to decision makers to take account of fuzzy effect while taking decision and by taking account of this; he/she earn more profit for the organization.

A future extension is to discuss model in more realistic situation by consider impreciseness in different inventory related cost and taking different form of demand pattern likes as time dependent, ram-type demand with inflation and permissible delay.

# A Two-Warehouse Inventory Model
# for Deteriorating Items with Stock Dependent Demand,
# Inflation and Genetic Algorithm

**Sanjey Kumar[1], Pallavi Agarwal[2] & Neeraj Kumar[3]\***

[1,3] Department of Mathematics, SRM University, Sonepat, Haryana

[2] Department of Computer Science & Engineering, SRM University, Sonepat, Haryana

*Corresponding Author, nkneerajkapil@gmail.com*

**Abstract**

In this article, the a two-warehouse inventory model deals with deteriorating items, with stock dependent demand rate  and model affected by inflation under the  pattern of  time value of money over a finite planning horizon.   Shortages are allowed and  partially backordered depending on the waiting time for the next replenishment. The purpose of this model is to minimize the total inventory cost by using Genetic algorithm. Also, a numerical example along with sensitivity analysis is given to explore the model numerically. Some observations are presented on the basis of sensitivity analysis.

**Keywords**: Two-warehouse, Genetic algorithm, partial backlogging, stock-dependent demand, Inflation, Deterioration, shortages;

## 1. Introduction

In the busy markets like super market, municipality market etc. the storage area of items is limited. When an attractive price discount for bulk purchase is available or the cost of procuring goods is higher than the other inventory related cost or demand of items is very high or there are some problems in frequent procurement, management decide to purchase a large amount of items at a time. These items cannot be accommodated in the existing storehouse (viz. the Own Warehouse, OW) located at busy market place. In the present senerio, suppliers proposed price discounts for gathering purchases or if the goods are seasonal, the retailers possibly will buy the superfluous goods that can be stored in own warehouse (OW). And rented warehouse (RW) is used as a store over the certain capacity $W_1$ of the own warehouse. Generally, the rented warehouse may have a costly superior unit holding cost than the own warehouse due to surplus cost of maintenance, material handling, etc. Hartely [1] was the original instigator to consider the impact of a two-warehouse model in inventory research and developed an inventory model with a RW storage principle. Sarma [2] developed a two-warehouse model for deteriorating items with an infinite replenishment rate and shortages. Sarma and Sastry [3] introduced a deterministic inventory model with an infinite production rate, permissible shortage and two levels of storage. Pakkala and Achary [4] considered a two-warehouse model for deteriorating items with finite replenishment rate and shortages. Lee and Ma [5] compared an optimal inventory policy for deteriorating items with two-warehouse and time-dependent demand. Yang [6] produced a two-warehouse inventory model with constant deteriorating items, constant demand rate and shortages under inflation. Yang [7] investigated the two-warehouse partial backlogging inventory models for deteriorating items under inflation. Kumar et.al [8] developed a Two-Warehouse inventory model without shortage for exponential demand rate and an optimum release rule. Kumar et al. [9] produced a Deterministic Two-warehouse Inventory Model for Deteriorating Items with Stock-dependent Demand and Shortages under the conditions of permissible delay. [10] Analyzed two-warehouse partial backlogging inventory models with three-parameter Weibull distribution deterioration under inflation. Sett et al. [11] introduced a two warehouse inventory model with increasing demand and time varying deterioration. Kumar et al. [12] developed Learning effect on an inventory model with two-level storage and partial backlogging under inflation. Yang and Chang [13] perused a two-warehouse partial backlogging inventory model for deteriorating items with permissible delay in payment under inflation. Guchhaita et al. [14] investigated a two storage inventory model of a deteriorating item with variable demand under partial credit period. Kumar and Singh [15] discussed Effect of Salvage Value on a Two-Warehouse Inventory Model for Deteriorating Items with Stock-Dependent Demand Rate and Partial Backlogging. Deterioration performed a most important contribution in lots of inventory systems. Normally, an inventory model understands with non-deteriorating items and instantaneous deteriorating items. Major part of goods undergo, waste or deterioration over time, examples being medicines, volatile liquids, blood banks, and so on. Therefore, waste or deterioration of physical goods in stock is a more realistic factor and there is a big need to consider inventory modeling. The primary effort to describe the optimal ordering policies for such items was prepared by Ghare and Schrader [16]. They presented an EOQ model for an exponentially decaying inventory. Philip [17] developed an inventory model with three parameter Weibull distribution rate without considering shortages. Deb and Chaudhari [18] derived inventory model with time-dependent deterioration rate. A meticulous assessment of deteriorating inventory literatures is given by Goyal and Giri [19]. Liao [20] studied an EOQ model with non- instantaneous receipt and exponential deteriorating item under two level trade credits. Chung [21] derived a complete proof on the solution procedure for non-instantaneous deteriorating items with permissible delay in payment. Chang et al. [22] framed optimal replenishment policies for non-instantaneous deteriorating items with stock-dependent demand. Dye [23] investigated the effect of .0preservation technology investment in a non-instantaneous deteriorating inventory model.

Due to high inflation and consequent sharp decline in the purchasing power of money in the developing countries like Brazil, Argentina, India, Bangladesh etc., the financial situation has been completely changed and so it is not possible to ignore the effect of inflation and time value of money any further. Following Buzacott [24] and Misra [25] have extended their approaches to different inventory models by considering the time value of money, different inflation rtes for the internal and external costs, finite replenishment, shortages, etc. Datta and Pal [26] considered the effects of inflation and time value of money of an inventory model with a linear time-dependent demand rate and shortages. Sarker and Pan [27] considered a finite replenishment model when the shortage is allowed. Chung [28] developed an algorithm with finite replenishment and infinite planning horizon. Tolgari et al. [29] studied an inventory model for imperfect items under inflationary conditions by considering inspection errors. Guria et al. [30] formulated an inventory policy for an item with inflation induced purchasing price, selling price and demand with immediate part payment. In the case of perishable product, the retailer may need to backlog demand to avoid costs due to deterioration. When the shortage occurs, some customers are willing to wait for back order and others would turn to buy from other sellers. Inventory model of deteriorating items with time proportional backlogging rate has been developed by Dye et al. [31]. Wang [32] studied shortages and partial backlogging of items. Recently, Kumar and Kumar [33] developed an inventory model with stock dependent demand rate for deterioration items.

In this study, we have developed an inventory model for non-instantaneous deteriorating items with stock-dependent under the impact of inflation with genetic algorithm. Shortages are allowed and partially backordered depending on the waiting time for the next replenishment. The main objective of this work is minimizing the total inventory cost and finding the new optimal interval and the optimal order quantity. The model shows the effect of the genetic algorithm due to changes in various parameters by taking suitable numerical examples and sensitivity analysis.

## 2. Notations and Assumptions

### 2.1 Notations

The following notations are used throughout this paper:

$A$   The ordering cost per order

$C_{hr}$    The ordering cost per item in RW

$C_{ho}$    The holding cost per item in OW, $C_{hr} > C_{ho}$

$C_2$    The deterioration cost per unit per unit cycle

$C_3$    The shortage cost for backlogged items per unit per unit cycle

$C_4$    The unit cost of lost sales per unit per cycle

$p$    The purchasing cost per unit

$s$    The selling price per unit, with $s > p$

$\mu_1$    The life time of the items in OW

$\mu_2$    The life time of the items in RW, $\mu_1 < \mu_2$

$\alpha$    The deterioration rate in OW, $0 \le \alpha < 1$

$\beta$    The deterioration rate in RW, $0 \le \beta < 1, \alpha > \beta$

$T$    The length of the order cycle (decision variable)

$H$    The planning horizon

$m$    The number of replenishment during planning horizon, m = H/T (decision variable)

$W_1$    The capacity of OW

$W_2$    The maximum inventory level in RW (decision variable)

$S$    The maximum inventory level per cycle (decision variable)

$BI$    The maximum amount of shortage demand to be backlogged (decision variable)

$Q$    The $2^{nd}, 3^{rd}, ...., m^{th}$ order size (Decision variable)

$r$    The discount rate represents the time value of money.

$f$    The inflation rate

$R$    The net discount rate of inflation i.e. R = r – f

$q_r(t)$    The inventory level in RW at time t

$q_o(t)$    The inventory level in OW at time t

$q_s(t)$    The negative inventory level at time t

$T_j$    The total time that elapsed upto and including the jth replenishment cycle (j = 1, 2, 3….)

$t_r$    Length of period during which inventory level reaches to zero in RW

$t_j$     The time at which the inventory level in OW in the j$^{th}$ replenishment cycle drop to zero

(j = 1, 2,...., m).

$T_j - t_j$     The time period when shortage occurs ( $j = 1, 2,..., m$ )

$TC_f$     The total cost for first replenishment cycle

$TC$     The total cost of the system over a finite planning horizon $H$

## 2.2 Assumptions:

To develop the mathematical model, the following assumptions are being made:

1. A single item is considered over the prescribed period of planning horizon.

2. There is no replacement or repair of deteriorated items takes place in a given cycle.

3. The lead time is zero.

4. Deterioration takes place after the life time of items. That is, during the fixed period, the Product has no deterioration. After that, it will deteriorate with constant rate.

5. The replenishment takes place at an infinite rate.

6. The effects of inflation and time value of money are considered.

7. The demand rate (a + bq$_r$(t)) is a stock dependent.

8. Shortages are allowed and partially backlogged. During the stock out period, the backlogging rate is variable and is dependent on the length of the waiting time for the next replenishment. So the Backlogging rate of negative inventory is, 1/ (1 + δ(T − t)), where δ is backlogging parameter $0 \le \delta \le 1$ and $(T − t)$ is waiting time ($t_j \le t \le T$) , ($j = 1, 2, …, m$). The remaining fraction (1− $B(t)$) is lost.

9. The OW has limited capacity of $W_1$ units and the RW has unlimited capacity. For economic reasons, the items of RW are consumed first and next the items of OW.

## 3. Formulation and solution of the model

Suppose with the purpose of the planning horizon H is divided into m equal parts of length T = H/m. Hence the reorder times over the planning horizon H are T$_j$ = jT (j=0, 1, 2, … , m). When the inventory is positive, demand rate is stock dependent, whereas for negative inventory, the demand is partially backlogged. The period for which there is no-shortage in each interval [jT, (j+1)T] is a fraction of the scheduling period T and is equal to kT (0 < k < 1). Shortages occur at time t$_j$ = (k+j-1)T,   (j = 1, 2, …, m) and are build up until time t = jT (j = 1, 2, …, m) before they are backordered. This model is demonstrated in Figure-1. The first replenishment lot size of S is replenished at T$_0$ = 0. W$_1$ units are kept in OW and the rest is stored in RW. The items of OW are consumed only after consuming the goods kept in RW. In the RW, during the time interval [0, $\mu_2$] , the inventory level   is decreasing only due to demand rate and the inventory level is dropping to zero owing to demand and deterioration during the time interval [$\mu_2$, $t_r$ ]. In OW, during the time interval [0, $\mu_1$], there is no change in the inventory level. However, the inventory $W_1$ decreases during [$\mu_1$,$t_r$] due to deterioration only, but during [ $t_r$, $t_1$], the inventory is depleted due to both demand and deterioration. By the time $t_1$, both warehouses are empty. Finally, during the interval [$t_1$, T], shortages occur and accumulate until t = $T_1$ before they are partially backlogged.

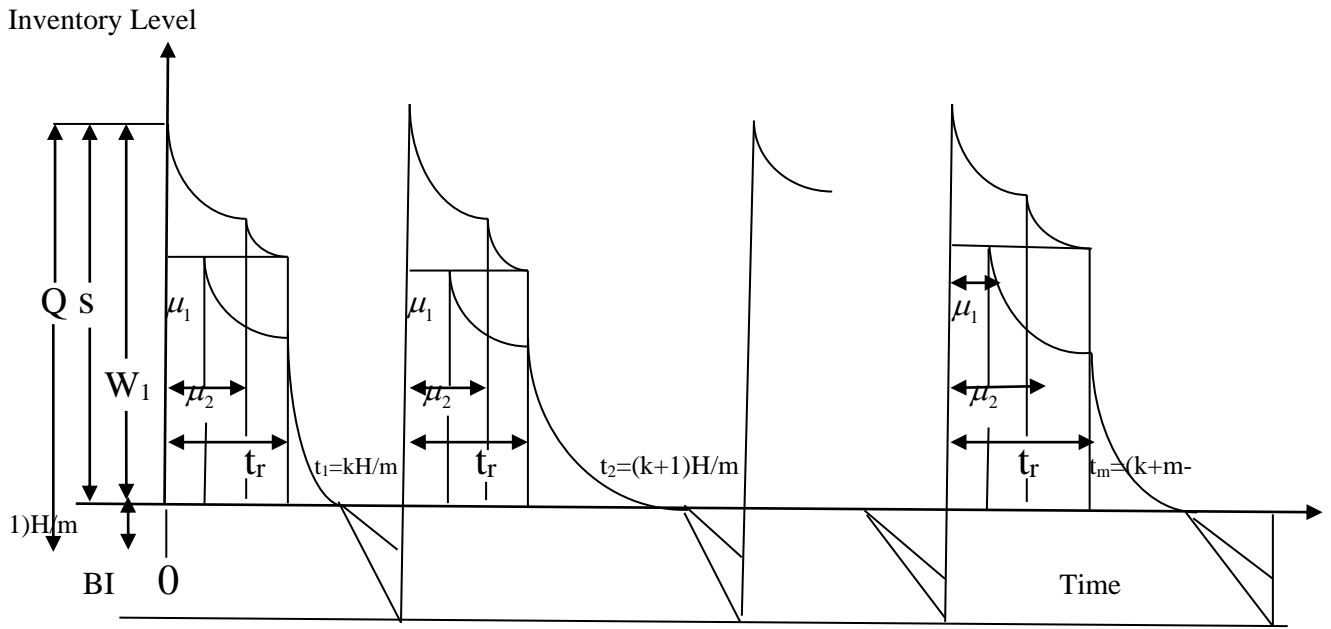

**Fig. 3.1. Graphical representation of the two warehouse inventory system**

Based on the above explanation during the time interval [0, $\mu_2$], the inventory level in RW is decreasing only due to demand rate and the differential equation representing the inventory status is given by

$$\frac{dq_r(t)}{dt} = -\left(a + bq_r(t)\right) \qquad 0 \leq t \leq \mu_2 \qquad (3.1)$$

With the condition $q_r(0) = W_2$, the solution of equation (3.1) is

$$q_r(t) = \frac{-a}{b} + \left(\frac{a}{b} + W_2\right)e^{-bt} \qquad (3.2)$$

In the second interval $[\mu_2, t_r]$ in RW, the inventory level decreases due to demand and deterioration. Thus, the differential equation below represents the inventory status

$$\frac{dq_r(t)}{dt} + \beta q_r(t) = -\big(a + bq_r(t)\big) \qquad \mu_2 \leq t \leq t_r \tag{3.3}$$

$$\frac{dq_r(t)}{dt} + \beta q_r(t) + bq_r(t) = -a \qquad \mu_2 \leq t \leq t_r \tag{3.4}$$

With the condition $q_r(t_r) = 0$, we get the solution of equation (3.4) is

$$q_r(t) = \frac{-a}{\beta+b} + \frac{a}{\beta+b} e^{(\beta+b)(t_r-t)} \qquad \mu_2 \leq t \leq t_r \tag{3.5}$$

Put $t = \mu_2$ in equations (3.2) and (3.5) we get the value of $W_2$ as

$$W_2 = \frac{a\beta}{b(\beta+b)}\left\{ e^{b\mu_2} - 1 - \frac{b}{\beta}\Big(1 - e^{(\beta+b)(t_r-\mu_2)+b\mu_2}\Big) \right\} \tag{3.6}$$

Putting the value of $W_2$ in equation (3.2) we get

$$q_r(t) = \frac{a}{b}\Big[-1 + e^{-bt}\Big] + \frac{a\beta}{b(\beta+b)}\left\{ e^{b\mu_2} - 1 - \frac{b}{\beta}\Big(1 - e^{(\beta+b)(t_r-\mu_2)+b\mu_2}\Big) \right\}, \qquad 0 \leq t \leq \mu_2 \tag{3.7}$$

In OW, during the interval $[0, \mu_1]$, there is no change in the inventory level and during $[\mu_1, t_r]$ the inventory $W_1$ only decreases due to deterioration.

Therefore the rate of change in the inventory is given by

$$\frac{dq_o(t)}{dt} = 0 \qquad 0 \leq t \leq \mu_1 \tag{3.8}$$

$$\frac{dq_o(t)}{dt} + \alpha q_o(t) = 0 \qquad \mu_1 \leq t \leq t_r \tag{3.9}$$

With the conditions $q_0(0) = W_1$ and $q_0(\mu_1) = W_1$, the solutions of equations (3.8) and (3.9) are

$$q_o(t) = W_1 \qquad 0 \leq t \leq \mu_1 \tag{3.10}$$

$$q_o(t) = W_1 e^{\alpha(\mu_1-t)} \qquad \mu_1 \leq t \leq t_r \tag{3.11}$$

In the interval $[t_r, t_1]$ in own warehouse, the inventory level decreases due to demand and deterioration. Thus, the differential equation is

$$\frac{dq_o(t)}{dt} + \alpha q_o(t) = -\left(a + bq_r(t)\right) \qquad\qquad t_r \le t \le t_1 \tag{3.12}$$

With the condition $q_0(t_1) = 0$, we get the solution of equation (3.12) is

$$q_o(t) = \frac{a}{\alpha + b}\left(e^{(\alpha+b)(t_1-t)} - 1\right) \qquad\qquad t_r \le t \le t_1 \tag{3.13}$$

Put $t = t_r$ in equations (3.11) and (3.13) we get,

$$t_r = \frac{\left(W_1 \alpha e^{\alpha\mu_1} + ae^{(\alpha+b)t_1}\right) \pm \sqrt{\left(W_1 \alpha e^{\alpha\mu_1} + ae^{(\alpha+b)t_1}\right)^2 + \frac{4a\alpha b}{\alpha+b}\left(W_1 e^{\alpha\mu_1} - \frac{a}{\alpha+b}e^{(\alpha+b)t_1} + \frac{a}{\alpha+b}\right)}}{2\left(\frac{-a\alpha b}{\alpha+b}\right)} \tag{3.14}$$

During the interval $[t_1, T]$, shortages occurred and the demand is partially backlogged. That is, the inventory level at time $t$ is governed by the following differential equation

$$\frac{dq_s(t)}{dt} = \frac{-a}{1 + \delta(T-t)} \qquad\qquad t_1 \le t \le T \tag{3.15}$$

With the condition $q_s(t_1) = 0$ the solution of equation (3.15) is

$$q_s(t) = a(t_1 - t)\left[1 - \delta T + \frac{\delta}{2}(t_1 + t)\right] \qquad\qquad t_1 \le t \le T \tag{3.16}$$

Therefore the maximum inventory level and maximum amount of shortage demand to be backlogged during the first replenishment cycle are

$$S = W_1 + \frac{a\beta}{b(\beta + b)}\left\{e^{b\mu_2} - 1 - \frac{b}{\beta}\left(1 - e^{(\beta+b)(t_r-\mu_2)+b\mu_2}\right)\right\} \tag{3.17}$$

$$BI = \frac{aH(1-k)}{2m^2}\left[2m - \delta H(1-k)\right] \tag{3.18}$$

There are m cycles during the planning horizon. Since, inventory is assumed to start and end at zero, an extra replenishment at $T_m = H$ is required to satisfy the backorders of the last cycle in the planning horizon. Therefore, there are m + 1 replenishments in the entire planning horizon H.

The first replenishment lot size is S.

The 2, 3, ..., $m^{th}$ replenishment order size is:- $\quad Q = S + BI$ (3.19)

The last or $(m + 1)^{th}$ replenishment lot size is BI.

Since replenishment in each cycle is done at the start of each cycle, the present value of ordering cost during the first cycle is $\quad OC=A$ (3.20)

The holding cost for the RW during the first replenishment cycle is

$$HC_r = C_{hr}\left[\int_0^{\mu_2} q_r(t)e^{-Rt}dt + \int_{\mu_2}^{t_r} q_r(t)e^{-Rt}dt\right]$$ (3.21)

$$= C_{hr}\left[-\frac{a}{b}\left\{\frac{1-e^{-R\mu_2}}{R}\right\} - \left(\frac{a}{b}+W_2\right)\left\{\frac{e^{-(R+b)\mu_2}-1}{R+b}\right\}\right]$$

$$+ C_{hr}\left[\frac{a}{\beta+b}\left\{e^{-Rt_r}\left(\frac{\beta+b}{R(\beta+b+R)}\right) + e^{-R\mu_2}\left(\frac{e^{(t_r-\mu_2)(\beta+b)}}{\beta+b+R}-\frac{1}{R}\right)\right\}\right]$$ (3.22)

The holding cost for the OW during the first replenishment

$$HC_o = C_{ho}\left[\int_0^{\mu_1} q_o(t)e^{-Rt}dt + \int_{\mu_1}^{t_r} q_o(t)e^{-Rt}dt + \int_{t_r}^{t_1} q_o(t)e^{-Rt}\right]$$

$$= C_{ho}\left[\frac{W_1}{R}\left(1-e^{-R\mu_1}\right) + \frac{W_1 e^{\alpha\mu_1}}{\alpha+R}\left(e^{-(\alpha+R)\mu_1}-e^{-(\alpha+R)t_r}\right)\right.$$

$$\left. + \frac{a}{\alpha+b}\left\{e^{-Rt_1}\left(\frac{\alpha+b}{R(\alpha+b+R)}\right) + e^{-Rt_r}\left(\frac{e^{(\alpha+b)(t_1-t_r)}}{\alpha+b+R}-\frac{1}{R}\right)\right\}\right]$$ (3.23)

The deteriorating cost for RW during first replenishment cycle is

$$DC_r = C_2\int_{\mu_2}^{t_r} \beta q_r(t)e^{-Rt}dt$$

$$= C_2 \frac{a\beta}{\beta+b} \left\{ e^{-Rt_r} \left( \frac{\beta+b}{R(\beta+b+R)} \right) + e^{-R\mu_2} \left( \frac{e^{(t_r-\mu_2)(\beta+b)}}{\beta+b+R} - \frac{1}{R} \right) \right\}$$

(3.24)

The deteriorating cost for OW during first replenishment cycle is

$$DC_o = C_2\alpha \left[ \int_{\mu_1}^{t_r} q_o(t)e^{-Rt}dt + \int_{t_r}^{t_1} q_o(t)e^{-Rt}dt \right]$$

$$= C_2\alpha \left[ \frac{W_1 e^{\alpha\mu_1}}{\alpha+R} \left( e^{-\mu_1(\alpha+R)} - e^{-(\alpha+R)t_r} \right) + \frac{a}{\alpha+b} e^{-Rt_1} \left( \frac{\alpha+b}{R(\alpha+b+R)} \right) \right.$$

$$\left. + \frac{a}{\alpha+b} e^{-Rt_r} \left( \frac{e^{(t_1-t_r)(\alpha+b)}}{\alpha+b+R} - \frac{1}{R} \right) \right]$$

(3.25)

Total shortage cost during the first replenishment cycle is

$$SC = -C_3 \int_{t_1}^{T} q_s(t)e^{-Rt}dt$$

$$= \frac{-C_3 a}{R} \left\{ e^{\frac{-RH}{m}} \left[ \frac{KH^2\delta}{m^2} \left( 1 - \frac{K}{2} - \frac{m}{\delta H} \right) + \frac{H}{m} \left( 1 - \frac{\delta H}{2m} \right) + \frac{1}{R} \left( 1 + \frac{\delta}{R} \right) \right] \right.$$

$$\left. + \frac{e^{\frac{-RKH}{m}}}{R} \left[ \frac{\delta H}{m} \left( 1 - K - \frac{m}{HR} \right) - 1 \right] \right\}$$

(3.26)

The lost sale cost during the first replenishment cycle is

$$LC = C_4 \int_{t_1}^{T} \left( 1 - \frac{1}{1+\delta(T-t)} \right) ae^{-Rt}dt$$

$$= \frac{-C_4 a\delta}{R^2} \left[ e^{\frac{-RKH}{m}} + e^{\frac{-RKH}{m}} \left( \frac{RH}{m}(1-K) - 1 \right) \right]$$

(3.27)

Replenishment is done at t = 0 and T. The present value of purchasing cost PC during the first replenishment cycle is:-

$$PC = pS + pe^{-RT}(BI)$$

$$= p\left\{W_1 + \frac{a\beta}{b(\beta+b)}\left\{e^{b\mu_2} - 1 - \frac{b}{\beta}\left(1 - e^{(\beta+b)(t_r-\mu_2)+b\mu_2}\right)\right\} + ae^{\frac{-RH}{m}}\frac{H}{M}(1-k)\left[1 + \frac{\delta H}{2m}(K-1)\right]\right\} \quad (3.28)$$

So, the total cost = Ordering cost + inventory holding cost in RW + inventory holding cost in OW + deterioration cost in RW + deterioration cost in OW + shortage cost + lost sales cost + purchasing cost.

$$TC_F = OC + HC_r + HC_0 + DC_r + DC_0 + SC + LC + PC$$

So, the present value of the total cost of the system over a finite planning horizon H is

$$TC(m,k) = \sum_{j=0}^{m-1} TC_f e^{-RjT} + Ae^{-RH} = TC_F\left(\frac{1-e^{-RH}}{1-e^{\frac{-RH}{m}}}\right) + Ae^{-RH} \quad (3.29)$$

Where T = H/m and TC$_F$ derived by substituting equations (3.21) to (3.28) in equation (3.29).

On simplification we get

$$TC(m,k) = Ae^{-RH} + G\left[A + C_{hr}\left\{\frac{a}{b}\left(\frac{e^{-R\mu_2}-1}{R}\right) + \left(\frac{a}{b}+W_2\right)\left(\frac{e^{-(R+b)\mu_2}-1}{R+b}\right)\right\}\right.$$

$$+\left(\frac{C_{hr}a}{\beta+b}+\frac{C_2 a\beta}{\beta+b}\right)\left\{\frac{e^{-Rt_r}(\beta+b)}{R(\beta+b+R)} + e^{-R\mu_2}\left(\frac{e^{(t_r-\mu_2)(\beta+b)}}{\beta+b+R}-\frac{1}{R}\right)\right\}$$

$$+C_{ho}\left[\frac{W_1}{R}\left(1-e^{-R\mu_1}\right) + \frac{W_1 e^{\alpha\mu_2}}{\alpha+R}\left(e^{-(\alpha+R)\mu_1}-e^{-(\alpha+R)t_r}\right)\right]$$

$$+\left(C_{ho}+C_2\alpha\right)\left[\frac{W_1 e^{\alpha\mu_2}}{\alpha+R}\left(e^{-\mu_1(\alpha+R)}-e^{-(\alpha+R)t_r}\right)\right]$$

$$+C_2\alpha\left[\frac{a}{R(\alpha+b+R)}e^{\frac{-RKH}{m}} + \frac{a}{\alpha+b}\left(\frac{e^{\left(\frac{KH}{m}-t_r\right)(\alpha+b)}}{\alpha+b+R}-\frac{1}{R}\right)e^{-Rt_r}\right]$$

$$-\frac{C_3 a}{R}\left[e^{\frac{-RH}{m}}\left\{\frac{KH^2\delta}{m^2}\left(1-\frac{K}{2}-\frac{m}{\delta H}\right)+\frac{H}{m}\left(1-\frac{\delta H}{2m}\right)+\frac{1}{R}\left(1+\frac{\delta}{R}\right)\right\}\right.$$

$$+\frac{e^{\frac{-RKH}{m}}}{R}\left\{\frac{\delta H}{m}\left(1-K-\frac{m}{HR}\right)-1\right\}\right] - \frac{C_4 a\delta}{R^2}\left[e^{\frac{-RKH}{m}} - e^{\frac{-RKH}{m}}\left(\frac{RH}{m}(1-K)-1\right)\right]$$

$$+P\left[W_1 + \frac{a\beta}{b(\beta+b)}\left\{e^{b\mu_2}-1-\frac{b}{\beta}\left(1-e^{(\beta+b)(t_r-\mu_2)+b\mu_2}\right)\right\}\right]$$

$$+ae^{\frac{-RH}{m}}\frac{H}{m}(1-K)\left(1-\frac{\delta H}{2m}(K-1)\right)\right] \qquad \text{Where } G = \left(\frac{1-e^{RH}}{1-e^{-RH/m}}\right) \qquad (3.30)$$

## 4. Solution Procedure

The present value of total cost TC (m, k) is a function of two variables m and k where m is a discrete variable and k is a continuous variable. For a given value of m, the necessary condition for TC(m, k) to be minimized is dTC (m, k ) / dk = 0 which gives

$$\frac{dTC(m,k)}{dk} = C_2\alpha\left[\frac{aH\ e^{\frac{kH}{m}(\alpha+b)}}{m(\alpha+b+R)} - \frac{aH\ e^{\frac{-RHk}{m}}}{m(\alpha+b+R)}\right] + C_3\left(\frac{aH}{Rm}e^{\frac{-RH}{m}}\right) + \frac{H^2\delta}{m^2}e^{\frac{-RHk}{m}}(k-1)$$

$$+\frac{H\delta}{mR}e^{\frac{-RHk}{m}} + \frac{H}{m}e^{\frac{-RHk}{m}} - C_4\frac{\delta a}{R^2}\left[\frac{R^2 H^2}{m^2}e^{\frac{-RHk}{m}}(1-k)-\frac{RH}{m}e^{\frac{-RHk}{m}}\right]$$

$$+P\left[(1-k)\frac{a\ H}{m}e^{\frac{-RH}{m}}\left\{1-\frac{H\delta}{m}(1+k)\right\}\right] \qquad (3.31)$$

$$\frac{d^2 TC(m,k)}{dk^2} = C_2 \alpha \left[ \frac{aH^2}{m^2} \left( \frac{R\, e^{\frac{-RHk}{m}}}{(\alpha+b+R)} + (\alpha+b) e^{\frac{kH}{m}(\alpha+b)} \right) \right] + \frac{H^3 \delta R}{m^3} e^{\frac{-RHk}{m}} (1-k)$$

$$-C_4 \frac{\delta a}{R^2} \left[ \frac{R^3 H^3}{m^3} e^{\frac{-RHk}{m}} (k-1) \right] + P \left[ \frac{a\delta H^2}{2m^2} \left( e^{\frac{-RH}{m}} - 1 \right) + \frac{aH}{m} e^{\frac{-RH}{m}} \left( \frac{\delta kH}{m} - 1 \right) \right] \quad (3.32)$$

Furthermore, the equation (32) shows that TC (m, k) is convex with respect to k. So, for a given positive integer m, the optimal value of k can be obtained from (31).

**Algorithm**

Step 1: Start with m = 1.

Step 2: Using (31) solve for k. Then substitute the solution obtained for (31) into (30) to compute the total inventory cost.

Step 3: Increase m by 1 and repeat step 2.

Step 4: Repeat step 2 and step 3 until TC (m, k) increases. The value of m which corresponds to the increase of TC for the first time is taken as the optimal value of m (denoted by m*) and the corresponding k (denoted by k*) is the optimal value of k.

Using the optimal solution procedure described above, we can find the optimal order quantity and maximum inventory levels to be

$$W_2 = \frac{a\beta}{b(\beta+b)} \left\{ e^{b\mu_2} - 1 - \frac{b}{\beta} \left( 1 - e^{(\beta+b)(t_r - \mu_2) + b\mu_2} \right) \right\}$$

$$S^* = W_1 + \frac{a\beta}{b(\beta+b)} \left\{ e^{b\mu_2} - 1 - \frac{b}{\beta} \left( 1 - e^{(\beta+b)(t_r - \mu_2) + b\mu_2} \right) \right\}$$

$$Q^* = W_1 + \frac{a\beta}{b(\beta+b)} \left\{ e^{b\mu_2} - 1 - \frac{b}{\beta} \left( 1 - e^{(\beta+b)(t_r - \mu_2) + b\mu_2} \right) \right\} + \frac{aH(1-k)}{2m^2} \left[ 2m - \delta H(1-k) \right]$$

Where $t_r = \dfrac{\left( W_1 \alpha e^{\alpha\mu_1} + ae^{(\alpha+b)t_1} \right) \pm \sqrt{\left( W_1 \alpha e^{\alpha\mu_1} + ae^{(\alpha+b)t_1} \right)^2 + \dfrac{4a\alpha b}{\alpha+b} \left( W_1 e^{\alpha\mu_1} - \dfrac{a}{\alpha+b} e^{(\alpha+b)t_1} + \dfrac{a}{\alpha+b} \right)}}{2\left( \dfrac{-a\alpha b}{\alpha+b} \right)}$

## 5. Numerical Examples

**Example 1**

Consider an inventory system with the following data: *D* = 100 units; $W_1$ = 50 units; *p* = \$4; *s* =

$15; $A$ = $150; $C_{hr}$ = $2; $C_{ho}$ = $1.2; $C_2$ = $1.5; $C_3$ = $5; $C_4$ = $10; $\alpha$ = 0.8; $\beta$ = 0.2; $\delta$ = 0.008; $\mu_1$ = 5/12 year; $\mu_2$ = 8/12 year; $R$ = 0.2; $H$ = 20 years.

## 6. Implementation of Genetic Algorithm

A genetic algorithm (GA) is a based on natural selection process to optimized tools that minimizes the total costs in supply chain management. It is a evolutionary computation method to solve a inventory problems. This is the more effective methods to find the optimized solution. The genetic algorithm uses three main types of rules at each step to create the next generation from the current population.

**The basic steps to find the optimized solution:**

**Step 1.** First one is ***Selection rules***, In this we select the individuals, called *parents* that contribute to the population at the next generation.

**Step 2.** Next one is ***Crossover rules***, *In this* we perform crossover operation between two parents to form children for the next generation.

**Step 3.** Last one is ***Mutation rules***, In mutation we apply some random changes to individual parents.

**We will perform these steps till we will not get our optimized solution.**

GA is not a method to find the exact solution of problem it only help to find the best or optimized solutions.

Here we are implementing Genetic Algorihtm in Table 1: Optimal total cost with respect to m and Table 2: Optimal total cost with respect to m when $\mu_{1=0 \text{ and }} \mu_{2=0}$

### Table 1: Optimal total cost with respect to m

| M | k(m) | $t_r$ | $t_1$ | T | Q | TC (m, |
|---|------|-------|-------|---|---|--------|
| 1 | 0.4485 | 8.9519 | 8.8552 | 19 | 3922 | 26496 |
| 2 | 0.4993 | 4.8536 | 4.6241 | 9.5 | 1482 | 11060 |
| **M** | **k(m)** | **$t_r$** | **$t_1$** | **T** | **Q** | **TC (m,** |
| 1 | 0.3582 | 1.9312 | 1.6842 | 2.76 | 452 | 5658 |
| 2 | 0.2995 | 0.8767 | 1.3582 | 2.653 | 409 | 5594 |

**Before Crossover:**

C1 ↓      C2↓

| M | k(m) | $t_r$ | $t_1$ | T | Q | TC (m, |
|---|---|---|---|---|---|---|
| 1 | 0.4485 | 8.9519 | 8.8552 | 19 | 3922 | 26496 |
| 2 | 0.4993 | 4.8536 | 4.6241 | 9.5 | 1482 | 11060 |

**After Crossover**

C1      C2

↓      ↓

| M | k(m) | $t_r$ | $t_1$ | T | Q | TC (m, |
|---|---|---|---|---|---|---|
| 1 | 0.4485 | 4.8536 | 4.6241 | 19 | 3922 | 26597 |
| 2 | 0.4993 | 8.9519 | 8.8552 | 9.5 | 1482 | 10090 |

**Mutation**

**Before mutation**    M1    M2

↓      ↓

| M | k(m) | $t_r$ | $t_1$ | T | Q | TC (m, |
|---|---|---|---|---|---|---|
| 1 | 0.4485 | 4.8536 | 4.6241 | 19 | 3922 | 26597 |

**After mutation**    M1    M2

↓      ↓

| M | k(m) | $t_r$ | $t_1$ | T | Q | TC (m, |
|---|---|---|---|---|---|---|
| 1 | 0.4485 | 4.6241 | 4.8536 | 19 | 3922 | 24384 |

**Before Crossover:**

C1 ↓          C2↓

| M | k(m) | $t_r$ | $t_1$ | T | Q | TC (m, k) |
|---|------|-------|-------|---|---|-----------|
| 1 | 0.3582 | 1.9312 | 1.6842 | 2.76 | 452 | 5658 |
| 2 | 0.2995 | 0.8767 | 1.3582 | 2.653 | 409 | 5594 |

**After Crossover**

C1          C2

↓          ↓

| M | k(m) | $t_r$ | $t_1$ | T | Q | TC (m, k) |
|---|------|-------|-------|---|---|-----------|
| 1 | 0.3582 | 0.8767 | 1.3582 | 2.76 | 452 | 5703 |
| 2 | 0.2995 | 1.9312 | 1.6842 | 2.653 | 409 | 5385 |

**Mutation**

**Before mutation**          M1          M2

↓          ↓

| M | k(m) | $t_r$ | $t_1$ | T | Q | TC (m, k) |
|---|------|-------|-------|---|---|-----------|
| 1 | 0.3582 | 0.8767 | 1.3582 | 2.76 | 452 | 5703 |

**After mutation**          M1          M2

↓          ↓

| M | k(m) | $t_r$ | $t_1$ | T | Q | TC (m, k) |
|---|------|-------|-------|---|---|-----------|
| 1 | 0.3582 | 1.3582 | 0.8767 | 2.76 | 452 | 5492 |

**Table 2: Optimal total cost with respect to *m when* $\mu_{1=0 \text{ and }} \mu_{2=0}$**

| M | k(m) | $t_r$ | $t_1$ | T | Q | TC (m, k) |
|---|------|-------|-------|---|---|-----------|
| 1 | 0.4170 | 8.4331 | 8.6052 | 19 | 4022 | 27496 |
| 2 | 0.4593 | 4.4536 | 4.6281 | 9 | 1532 | 11865 |
| M | k(m) | $t_r$ | $t_1$ | T | Q | TC (m, k) |
| 1 | 0.2782 | 0.8312 | 2.0042 | 1.763 | 612 | 6758 |
| 2 | 0.2098 | 0.5167 | 1.3249 | 1.953 | 456 | 6690 |

**Before Crossover:**

C1 ↓          C2 ↓

| M | k(m) | $t_r$ | $t_1$ | T | Q | TC (m, k) |
|---|------|-------|-------|---|---|-----------|
| 1 | 0.4170 | 8.4331 | 8.6052 | 19 | 4022 | 27496 |
| 2 | 0.4593 | 4.4536 | 4.6281 | 9 | 1532 | 11865 |

**After Crossover**        C1          C2

↓          ↓

| M | k(m) | $t_r$ | $t_1$ | T | Q | TC (m, k) |
|---|------|-------|-------|---|---|-----------|
| 1 | 0.4170 | 4.4536 | 4.6281 | 19 | 4022 | 26585 |
| 2 | 0.4593 | 8.4331 | 8.6052 | 9 | 1532 | 10336 |

**Mutation**

**Before mutation**        M1          M2

↓          ↓

| M | k(m) | $t_r$ | $t_1$ | T | Q | TC (m, k) |
|---|------|-------|-------|---|---|-----------|
| 1 | 0.4170 | 4.4536 | 4.6281 | 19 | 4022 | 26585 |

**After mutation**   **M1**    **M2**

| M | k(m) | t_r | t_l | T | Q | TC (m, k) |
|---|------|-----|-----|---|---|-----------|
| 1 | 0.4170 | 4.6281 | 4.4536 | 19 | 4022 | 24348 |

**Before Crossover:**

**C1**    **C2**

| M | k(m) | t_r | t_l | T | Q | TC (m, k) |
|---|------|-----|-----|---|---|-----------|
| 1 | 0.2782 | 0.8312 | 2.0042 | 1.763 | 612 | 6758 |
| 2 | 0.2098 | 0.5167 | 1.3249 | 1.953 | 456 | 6690 |

**After Crossover**

**C1**    **C2**

| M | k(m) | t_r | t_l | T | Q | TC (m, k) |
|---|------|-----|-----|---|---|-----------|
| 1 | 0.2782 | 0.5167 | 1.3249 | 1.763 | 612 | 6597 |
| 2 | 0.2098 | 0.8312 | 2.0042 | 1.953 | 456 | 6385 |

**Mutation**

**Before mutation**   **M1**    **M2**

| M | k(m) | t_r | t_l | T | Q | TC (m, k) |
|---|------|-----|-----|---|---|-----------|
| 1 | 0.2782 | 0.5167 | 1.3249 | 1.763 | 612 | 6597 |

**After mutation**      **M1**        **M2**

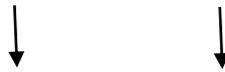

| M | k(m) | $t_r$ | $t_1$ | T | Q | TC (m, k) |
|---|------|-------|-------|---|---|-----------|
| 1 | 0.2782 | 1.3249 | 0.5167 | 1.763 | 612 | 6376 |

It is clearly visible the changes in total cost after the mutation methods and before crossover and it are very effective cost.

**Based on our numerical results, we obtain the following decision-making phenomena:**

(i) The total cost is increasing while the time horizon H is increasing, the order quantity is decreasing. Furthermore it is observed that when the time horizon increases, the number of replenishment also increases. Therefore the ordering cost increases. In order to minimize the total cost, the retailer should decrease the time horizon.

(ii) When the total cost of the retailer and the order quantity are also increasing then the deterioration rate β in RW increasing. When the deterioration rate will increase, the total cost of the retailer will increase.

3. If backlogging parameter δ is increased then the total cost and the order quantity will be decreased. But there is no change in the number of replenishment. If the backlogging parameter increases, then the ordering quantity will decrease.

4. When the inflation rate increases, the number of replenishment also increases. then the net discount rate of inflation R is increasing, the optimal cost is decreasing and the order quantity is also decreasing.

5. When the purchasing price p is increasing, the total optimal cost and the order quantity are highly increasing. When the purchase price increases, the number of replenishment decreases. Also, the increasing of purchasing price will increase the total cost of the retailer.

**7. Conclusion:**

In the present model, we have considered two warehouse inventory models depending on the waiting time for the next replenishment. And shortages and partially backlogging are allowed. In this model, demand rate considered as stock dependent with inflation and model affected by Genetic Algorithm. An algorithm is designed to find the optimum solution of the proposed model. The model shows that the minimum time horizon will minimize the total cost of the retailer. Furthermore, sensitivity analysis is carried out with respect to the key parameters and helpful decision-making insights are obtained. The graphical illustrations are also given to analyze the efficiency of the model clearly. The proposed model incorporates some practical features that are likely to be linked with some kinds of inventory. additionally this model can be adopted in the inventory control of retail business such as food industries, seasonable cloths, domestic goods, automobile, electronic components etc. The proposed model can be extended in several ways. Like incorporate some more realistic features, such as quantity discount, multi item, trade credit strategy, etc., when the net discount rate of inflation and the backlogging rate are increased then the optimal total cost will be decreased. Also, sensitivity analysis of the model with respect to numerous system parameters has been carried out a number of decision-making.

The main aim of the present book is to suggest some improved estimators using auxiliary and attribute information in case of simple random sampling and stratified random sampling and some inventory models related to capacity constraints.

This volume is a collection of six papers, written by five co-authors (listed in the order of the papers): Dr. Rajesh Singh, Dr. Sachin Malik, Dr. Florentin Smarandache, Dr. Neeraj Kumar, Mr. Sanjey Kumar & Pallavi Agarwal.

In the first chapter authors suggest an estimator using two auxiliary variables in stratified random sampling for estimating population mean. In second chapter they proposed a family of estimators for estimating population means using known value of some population parameters. In Chapter third an almost unbiased estimator using known value of some population parameter(s) with known population proportion of an auxiliary variable has been used. In Chapter four authors investigates a fuzzy economic order quantity model for two storage facility. The demand, holding cost, ordering cost, storage capacity of the own - warehouse are taken as a trapezoidal fuzzy numbers and in Chapter five a two-warehouse inventory model deals with deteriorating items, with stock dependent demand rate and model affected by inflation under the pattern of time value of money over a finite planning horizon. Shortages are allowed and partially backordered depending on the waiting time for the next replenishment. The purpose of this model is to minimize the total inventory cost by using Genetic algorithm.

This book will be helpful for the researchers and students who are working in the field of sampling techniques and inventory control.

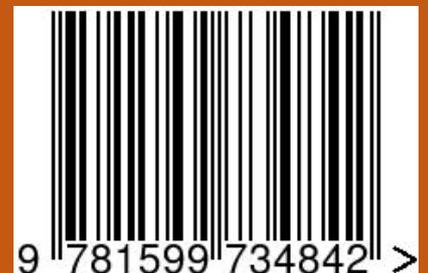